# The quantization conjecture revisited

By Constantin Teleman

**Abstract**

A strong version of the *quantization conjecture* of Guillemin and Sternberg is proved. For a reductive group action on a smooth, compact, polarized variety $(X, \mathcal{L})$, the cohomologies of $\mathcal{L}$ over the GIT quotient $X/\!/G$ equal the invariant part of the cohomologies over $X$. This generalizes the theorem of [GS] on global sections, and strengthens its subsequent extensions ([JK], [M]) to Riemann-Roch numbers. Remarkable by-products are the invariance of cohomology of vector bundles over $X/\!/G$ under a small change in the defining polarization or under shift desingularization, as well as a new proof of Boutot's theorem. Also studied are equivariant holomorphic forms and the equivariant Hodge-to-de Rham spectral sequences for $X$ and its strata, whose collapse is shown. One application is a new proof of the Borel-Weil-Bott theorem of [T1] for the moduli stack of $G$-bundles over a curve, and of analogous statements for the moduli stacks and spaces of bundles with parabolic structures. Collapse of the Hodge-to-de Rham sequences for these stacks is also shown.

**Contents**

# 0. Introduction

Associated to a linear action of a reductive group $G$ on a projectively embedded complex manifold $X$ there is a $G$-invariant stratification by locally closed, smooth subvarieties. The open stratum is the semistable locus $X^{\mathrm{ss}}$; the other, unstable strata were described, algebraically and symplectically, by Kirwan [K1] and Ness [N], based on algebraic work of Kempf [Ke] and Hesselink [He], and on topological ideas of Atiyah and Bott [AB]. The algebraic description uses the projective embedding, but the outcome depends only on the (equivariant) polarization class of $\mathcal{O}(1)$. Geometrically, a line bundle $\mathcal{L}$ with a positive hermitian metric, invariant under a compact form $K$ of $G$, defines a Kähler structure on $X$, and we are looking at the Morse stratification for the square-norm $\|\mu\|^2$ of the moment map for $G$.

Each space $H^*_S(X;\mathcal{L})$ of coherent sheaf cohomology with support on a stratum $S$ carries a natural action of $G$. The main observation of this paper[1] is the vanishing of its invariant part, for unstable $S$. A descending ordering $\{S(m)\}_{m\in\mathbb{N}}$ of the strata (where each union of $S(k)$, with $k\geq m$, is closed) leads to a *Cousin-Grothendieck spectral sequence*

$$E_1^{m,n} = H^{m+n}_{S(m)}(X;\mathcal{L}) \Rightarrow H^{m+n}(X;\mathcal{L}) \ ,$$

and the vanishing of invariants for positive $m$ implies that $H^*(X;\mathcal{L})^G = H^*(X^{\mathrm{ss}};\mathcal{L})^G$. The latter equals $H^*(X/\!/G;\mathcal{L})$, and we obtain a strong form of Guillemin and Sternberg's "quantization commutes with reduction" conjecture, which, based on their result for $H^0$, predicted the equality of the two holomorphic Euler characteristics.

That form of the conjecture, extended to the Spin$^c$-Dirac index on compact symplectic manifolds, has been proved in various degrees of generality: [V] for abelian groups, and [M], [JK] for smooth or orbifold quotients (see also [Sj2] for a survey); localization formulae were used to compute the two indices. A more conceptual (if analytically more involved) proof was given by Y. Tian and W. Zhang [TZ], using a Wittenesque deformation of the Dirac operator. After my paper was first circulated, their treatment of smooth Kähler quotients was refined by M. Braverman [Br] to give dimensional equality of the respective Dolbeault cohomologies. However, one does not quite get canonical isomorphisms this way. Also, the difficulties created by truly singular quotients are especially acute in the analytic treatment. In the symplectic approach of Meinrenken-Sjamaar, [MS], singular quotients are replaced by certain *partial* and *shift desingularizations* (Zhang [Z] has recently extended Braverman's argument along these lines).

---

[1] A special case of which was already noted by Ramadas [R]

Regularity of the quotient is not relevant to my argument, and the result holds even when $X$ itself has rational singularities (2.13). Remarkable *rigidity theorems* (5.5), (5.6) follow, asserting the invariance of the cohomology of holomorphic vector bundles under perturbation of the GIT quotient (by a small change of polarization or shift desingularization). This justifies *a posteriori* the avoidance of seriously singular quotients, and recovers an important theorem Boutot's (5.7) on rational singularities of quotients.[2]

The invariant part of coherent sheaf cohomology is an instance of *equivariant sheaf cohomology* $H_G^*$, definable for the action of any linear algebraic group (cf. Appendix). In general, it also involves higher group cohomology. Alternatively, we are discussing cohomologies over the *quotient stack* $X_G$ of $X$ by $G$, with supports on the locally closed substacks $S_G$. This point of view allows one to change groups and spaces as needed, via the following version of *Shapiro's lemma*: If $G$ is a subgroup of $G'$ and $X'$ is the *induced space* $G' \times^G X$, the quotient stacks $X_G$ and $X'_{G'}$ are equivalent, in the sense that $G$-equivariant computations on $X$ are equivalent to $G'$-equivariant ones on $X'$.

If the action of $G$ on $X^{\mathrm{ss}}$ is free, the quotient stack $X_G^{\mathrm{ss}}$ is the quotient variety; this is also the GIT *quotient* $X//G$ of $X$ by $G$. In the *stable case*, $X^{\mathrm{ss}}$ only carries finite isotropy groups, and $X_G^{\mathrm{ss}}$ is a *Deligne-Mumford stack* (A.8). It is closely related to the GIT quotient $X^{\mathrm{ss}}/G$, a compact Kähler space arising from the stack by ignoring the isotropies. The stable case is nearly as good as the free one, although there arises the delicate question whether equivariant line bundles over $X^{\mathrm{ss}}$ are pulled back from (*descend to*) $X//G$: line bundles over a DM stack may define *fractional* line bundles on its quotient space. When $\mathcal{L}$ does not descend, the statement of the theorem involves the *invariant direct image* sheaf $q_*^G \mathcal{L}$.

The cohomologies of $\mathcal{L}$ over the GIT and stack quotients also agree in the presence of positive-dimensional stabilizers. *Negative* powers of $\mathcal{L}$ are handled by Serre duality, but care is then needed with singular quotients, when the dualizing sheaves on $X$ and $X//G$ may not relate as naively expected (3.7). The naive form of the quantization conjecture requires a mild restriction (3.6); cf. also (3.8). The referee rightly noted that, according to [MS, §2.14], one gets a better result for the "$2\rho$–shifted" quotient, and this has been included in an addendum to Section 3.

For a torus, the invariants in $H_S^{>0}(X; \mathcal{K} \otimes \mathcal{L})$ vanish as well, and Kodaira's theorem on $X$ forces the vanishing of higher invariant cohomology of $K \otimes \mathcal{L}$ over $X^{\mathrm{ss}}$. Occasional failure of this for non-abelian $G$ is related to the difficulties caused by negative line bundles. In remedy, some vanishing conditions over singular quotients are given in Section 6. Sometimes, vanishing follows directly from the theorem of Grauert and Riemenschneider [GR]: thus, [KN] handled the moduli space of semistable principal $G$-bundles over a Riemann surface, but parabolic structures had only been treated for $SL_2$, in [MR] (by Frobenius splitting). They are now addressed uniformly in Section 9.

---

[2] In the other direction, Broer and Sjamaar [Sj1, Thm. 2.23] obtain a special case of the quantization conjecture from Boutot's theorem. A relative version of their argument recovers the rigidity theorems; cf. Section 5.

Section 7 introduces holomorphic forms. In support of the stack philosophy, the "quantization theorems" (7.1) and (7.3) apply to Kähler differentials *of the stack*, not to those of the variety. (We introduce them as equivariant differentials.) Twisting by $\mathcal{L}$ kills unstably supported cohomologies, but the more interesting result concerns untwisted forms. Collapse of the Cousin sequence is now related to the collapse at $E_1$ of the equivariant Hodge-to-de Rham spectral sequence. The latter is quite clear when $X$ is proper, but surprisingly, a weaker "KN completeness" condition suffices (§1). In the stable case, collapse for $X^s = X^{ss}$ is equivalent to collapse for the DM quotient stack, and well-known [St]; but the general case seems new.

Sections 8 and 9 apply these results to the stack of holomorphic principal bundles over a Riemann surface (enhanced with parabolic structures). Their open substacks of finite type can be realized as quotients of smooth quasi-projective varieties (the method goes back to Gieseker [G]). The cohomologies of suitable line bundles (conveniently mislabeled "positive") on the stacks equal those on the semistable moduli spaces; in fact, all higher cohomologies vanish. (Equality of the spaces of sections was already known from [BL] and [KNR].) On the stack side, this recovers a key part of the "Borel-Weil-Bott" theorem in [T1]. On the space side, it extends the vanishing theorem of [KN] to moduli of parabolic bundles.

Knowledgeable readers will notice the absence of new ideas in this paper. Indeed, both the question and the answer have been around for fifteen years, as several people came within a whisker of noticing ([R], [W]). The overlap in time with the independent proof [TZ] + [Br] seems entirely fortuitous: the latter draws on a completely different circle of ideas.

*Notes.* (i) This approach to the quantization conjecture was proposed in [T1], in connection with moduli of $G$-bundles. However, as T. R. Ramadas kindly pointed out, much of my argument in Section 2 had already appeared in [R], where vanishing of the group $H_S^{\mathrm{codim}\, S}(X; \mathcal{L})^G$ was shown. (Most ingredients for the full result — vanishing in all degrees — were in place, but only a weaker conclusion was drawn.)

(ii) In the symplectic case, there should be defined an "equivariant index with supports" of the pre-quantum line bundle, additive for the $\|\mu\|^2$-Morse stratification, taking values in the $\mathbb{Z}$-dual of the representation ring of $G$ (infinite sums of representations, with finite multiplicities). The invariant part of this index, with supports on unstable strata, should vanish. This approach to the symplectic theorem would allow for a relative version (cf. 5.2 and 5.3.i); it is not clear how the existing arguments can address that.

*Acknowledgements.* I am indebted to Shrawan Kumar for conversations on $G$-bundles. Hodge-theoretic exchanges with Carlos Simpson were also greatly helpful. I thank Y. Tian for acquainting me with his most recent work (joint with W. Zhang). Discussions and correspondence with Y. Hu, J. Kollár, J. Li, K. Liu, T. R. Ramadas, N. Shepherd-Barron, R. Sjamaar, C. Sorger and C. Woodward are gratefully acknowledged. The work was supported by Saint John's College, Cambridge, and by an NSF postdoctoral fellowship.

## 1. The strata of a projective $G$-variety

*The stratification defined by Kirwan and Ness.* Let $X, \mathcal{L}, G$ be as in the abstract. For simplicity, $X$ will always be irreducible. The open stratum $X^{ss} \subset X$ is the complement of the vanishing locus of $G$-invariant regular sections of large powers of $\mathcal{L}$. The other strata depend on the choice of a rational, invariant inner product in $\mathfrak{g}$. Fix (for convenience only) a Cartan subgroup $H$ of $G$ and a dominant Weyl chamber in $\mathfrak{h}_{\mathbb{R}}^t$. Any subtorus $T \subseteq H$ acts on the fibers $\mathcal{L}$, over each component $C$ of its fixed-point set $X^T$, by a character, which defines a rational weight $\beta$ of $\mathfrak{h}$, using the inner product. If $\beta \neq 0$, call $\mathbb{T} : \mathbb{C}^\times \to H$ the corresponding 1-parameter subgroup $Z$, that component of $X^{\mathbb{T}}$ containing $C$, and $L$ the commutant of $\mathbb{T}$ in $G$. Divide the natural action of $L$ on $(Z, \mathcal{L})$ fiberwise by $\beta$. (Raising $\mathcal{L}$ to some power ensures integrality of $\beta$ and $\mathbb{T}$, without affecting the construction to follow.)

The unstable strata are indexed by those $Z$ with dominant $\beta$ for which the semistable locus $Z^\circ \subseteq Z$ of the divided $L$-action on $\mathcal{L}$ is not empty. For such a $Z$, call $Y$ the set of points in $X$ flowing to $Z$ under $\mathbb{T}$, as $t \to \infty$ in $\mathbb{C}^\times$, and $Y^\circ$ the open subset flowing to $Z^\circ$.

(1.1) Properties (i)–(iv) below were proved in [K1]; (v) is from [MFK, Prop. 1.10].

(i) $Y$ is a fiber bundle over $Z$, with affine spaces as fibers, under the morphism $\varphi$ defined by the limiting value of the $\mathbb{T}$-flow.

(ii) $Y$ is stabilized by the parabolic subgroup $P \subset G$ whose nilpotent Lie algebra radical $\mathfrak{u}$ is spanned by the negative $\mathbb{T}$-eigenspaces in $\mathfrak{g}$.

(iii) The $G$-orbit $S$ of $Y^\circ$ is isomorphic to $G \times^P Y^\circ$. Under $\varphi$, it fibers in affine spaces over $G \times^P Z^\circ$, if we let $P$ act on $Z^\circ$ via its reductive quotient $L$.

(iv) The various $S$, together with $X^\circ = X^{ss}$, smoothly stratify[3] $X$.

(v) $Z^\circ$ has a projective, good (cf. 3.1) quotient under $L$; $X^{ss}$ has a good projective quotient under $G$.

Kirwan also showed that the $S$ are the Morse strata for $\|\mu\|^2$, in a $K$-invariant Kähler structure representing $c_1(\mathcal{L})$, and that $X^{ss}/G$ agreed with the "symplectic quotient" $\mu^{-1}(0)/K$ of $X$ by $K$.

KN *stratifications*. If $X$ is quasi-projective, statements (iv) and (v) can fail; for instance, deletion of part of some $Z$ leaves part of $S$ unaccounted for by $\varphi$. On the other hand, the prescription for the $\mathbb{T}$ can be relaxed. Hence, the following terminology seems useful. Consider a selection of one-parameter subgroups $\mathbb{T}$ in $H$, together with an $L$-invariant open $Z^\circ$ in the fixed-point set of each. Assume that the $G$-orbits $S$ of the sets $Y^\circ \subset X$ of points flowing to $Z^\circ$ under $\mathbb{T}$ satisfy (iii) and that, together with their complement $X^\circ$, assumed open, they stratify $X$. We call this a *Kirwan-Ness* (KN) *stratification*. Note that (i) and (ii) are automatic. If (v) holds for all $Z^\circ$, and also for $X^\circ$ with $L = G$, we say that the KN stratification is *complete*. A vector bundle is *adapted* (*strictly adapted*) to the stratification if its $\mathbb{T}$-weights on the fiber over the $Z$'s are nonnegative (positive). The relevance of these conditions to the "quantization theorem" was already identified in [TZ, 4.2].

*Example.* An open union of KN strata in a projective $X$ inherits a complete stratification; $\mathcal{O}$ is adapted, $\mathcal{L}$ strictly so. Another (analytic) example will be the Atiyah-Bott stratification; see (8.8).

*Change of polarization.* The stratification depends only on the *equivariant polarization* defined by $\mathcal{L}$, the line through its Chern class in $H_G^2(X; \mathbb{Q})$: this follows from (i)–(iv), given the result for semistable strata [MFK, 1.20]. (Morse theory also makes this clear, since line bundles in the same polarization class carry invariant metrics with the same curvature.) The effect of perturbing the polarization was first described in [DH, 3.3.15], although the idea is implicit in [K2, §3], while semistable strata were already discussed in [S1, §5]. The formulation below covers all the cases we need. Consider a projective $G$-morphism $\pi : X' \to X$, with relatively ample $G$-line bundle $\mathcal{M}$. For small positive $\varepsilon \in \mathbb{Q}$, $\mathcal{L}_\varepsilon : \pi^*\mathcal{L} + \varepsilon \cdot \mathcal{M}$ is an ample fractional $G$-line bundle.

(1.2) REFINEMENT LEMMA. *The $\mathcal{L}_\varepsilon$-stratification on $X'$ is independent of the small $\varepsilon > 0$, and refines the pull-back of the $\mathcal{L}$-stratification on $X$. Further, $\pi^*\mathcal{L}$ is adapted to this refined stratification.*

---

[3] Some authors use the term *decomposition* rather than *stratification*; see [K1, Def. 2.11] or [FM] for the reason.

(1.3) *Examples.* Three cases are especially important, and will be taken up in Section 5:

(i) $X' = X$, $\mathcal{M}$ is any $G$-bundle; this is the "change of polarization" studied in [DH] and [Th].

(ii) $X$ is singular, $X'$ is an equivariant desingularization. (The proof does not use smoothness of $X$; for singular varieties, stratification and moment map can be defined using the projective embedding.) The blow-up along an invariant subvariety was studied in [K2]; a natural choice for $\mathcal{M}$ is minus the exceptional divisor.

(iii) $X' = X \times F_\lambda$, $\mathcal{M} = \mathcal{O}(\lambda)$ over the flag variety $F_\lambda$ of $G$ corresponding to a dominant weight $\lambda$. This construction, going back at least to Seshadri [S1], is sometimes called *shift desingularization* of $X/\!/G$. If $G$ acts freely on $X^{ss}$, $X'/\!/G$ is an $F_\lambda$-bundle over $X/\!/G$.

(1.4) *Remark.* In all cases, the semistable stratum in $X'$ could be empty, even if $X^{ss}$ was not so; but this cannot happen if $X$ contains *stable* points [S1].

*Proof.* The stratification on $X'$ changes for finitely many values of $\varepsilon$ ([DH, 1.3.9]). Recall the idea: there are only finitely many possible $Z$'s and finitely many $Z^{ss}$ for each ([DH, 3.3.3], or argue inductively); and, given $Z^{ss}$, there are finitely many possible $Y^\circ$, as $\mathbb{T}$ varies ([DH, 1.3.8]). ($Y$ and $Z$ can only change upon vanishing of a $\mathbb{T}$-eigenvalue in the normal bundle to some fixed-point set of $H$.)

By [K1, 3.2], the weight $\beta$ associated to the stratum of $x \in X$ lies in the unique coadjoint orbit closest to zero in $\mu(\overline{Gx})$ (where $\overline{Gx}$ is the $x$-orbit closure, a compact set). If $x'$, $y' \in X'$ stay in the same stratum when $\varepsilon \to 0$, then $\beta_\varepsilon(x') = \beta_\varepsilon(y')$, so equality persists at $\varepsilon = 0$. Also, the two points flow to the same fixed-point component under $\mathbb{T}_\varepsilon$, and then also under $\mathbb{T}$; so they must lie over the same stratum in $X$.

Since $\mathcal{L}_\varepsilon$ is strictly adapted to the stratification and the $\beta_\varepsilon$ are $\varepsilon$-close to the $\beta$, $\pi^*\mathcal{L}$ is adapted. Moreover, over unstable strata of $X$, $\beta \neq 0$, so the (upstairs) $\mathbb{T}_\varepsilon$-weights on $\pi^*\mathcal{L}$ must be positive. □

## 2. Vanishing of cohomology with supports

Cohomology with supports on an unstable stratum $S$, of codimension $c$, can be rewritten as cohomology over $Z^\circ$. For a vector bundle $\mathcal{V}$ defined near $S$ in $X$, we have

(2.1) $$H_S^{*+c}(X; \mathcal{V}) = H^*(X; \mathcal{H}_S^c(\mathcal{V})) = H^*(S; \mathcal{R}_S\mathcal{V})$$

where the sheaf $\mathcal{R}_S\mathcal{V}$ over $S$ of $\mathcal{V}$-valued residues along $S$ in $X$ pushes forward to the local cohomology sheaf $\mathcal{H}_S^c(\mathcal{V}) = \lim_{n\to} \mathcal{E}xt_{\mathcal{O}}^c(\mathcal{O}/\mathcal{J}^n; \mathcal{V})$, on an open set $U \subset X$ in which $S$ is closed, with ideal sheaf $\mathcal{J}$. While $\mathcal{R}_S\mathcal{V}$ has no natural $\mathcal{O}_S$-module structure, it is increasingly filtered, and the composition series quotients are vector bundles over $S$. If $\mathcal{V}$ carries a $G$-action, so does $\mathcal{R}_S\mathcal{V}$, and the filtration is equivariant. It is for $\mathrm{Gr}\mathcal{R}_S\mathcal{V}$ that we shall prove the vanishing of invariant cohomology; the same then follows for $\mathcal{R}_S\mathcal{V}$.

Shapiro's lemma (A.5) equates the $G$-invariant part of the $E_1$ term in this spectral sequence with the $P$-equivariant cohomology $H_P^*(Y^\circ; \mathrm{Gr}\mathcal{R}_S\mathcal{V})$ (where the $\mathcal{O}$-module restriction of $\mathrm{Gr}\mathcal{R}_S\mathcal{V}$ to $Y^\circ$ is implied). Under $\varphi$, this equals $H_P^*(Z^\circ; \varphi_*\mathrm{Gr}\mathcal{R}_S\mathcal{V})$, which is the abutment of a spectral sequence (A.4)

$$(2.2) \qquad E_2^{r,s} = H_P^r\left(H^s\left(Z^\circ; \varphi_*\mathrm{Gr}\mathcal{R}_S\mathcal{V}\right)\right) \Rightarrow H_P^{r+s}\left(Z^\circ; \varphi_*\mathrm{Gr}\mathcal{R}_S\mathcal{V}\right) \ .$$

For a $P$-representation $V$, there are natural isomorphisms

$$H_P^*(V) = H^*(\mathfrak{p}, \mathfrak{l}; V) = H^*(\mathfrak{u}; V)^L$$

($L$ acts naturally on $V$, and by conjugation on $\mathfrak{u}$). Thus, we can rewrite, in (2.2),

$$(2.3) \qquad E_2^{r,s} = H^r\left(\mathfrak{u}; H^s\left(Z^\circ; \varphi_*\mathrm{Gr}\mathcal{R}_S\mathcal{V}\right)\right)^L \ ,$$

which is resolved by the $L$-invariant part of the Lie algebra Koszul complex for $\mathfrak{u}$,

$$(2.4) \qquad \left(H^s\left(Z^\circ; \varphi_*\mathrm{Gr}\mathcal{R}_S\mathcal{V}\right) \otimes \Lambda^r(\mathfrak{u})^t\right)^L \ .$$

When $\mathbb{T} \subset L$, $\varphi_*\mathrm{Gr}\mathcal{R}_S\mathcal{V}$ is $\mathbb{T}$-isomorphic to

$$(2.5) \qquad \mathcal{V} \otimes \mathrm{Sym}^\bullet(T_Z Y)^t \otimes \mathrm{Sym}^\bullet(T_S X) \otimes \det(T_S X) \ ;$$

the first two factors form the Gr of the fiberwise sections of $\mathcal{V}$ along $\varphi$ (filtered by the order of vanishing along $Z^\circ$), while the last two are $\mathrm{Gr}\mathcal{R}_S\mathcal{O}$. There follows the key observation of the paper.

(2.6) PROPOSITION. (a) $H_S^*(X; \mathcal{L})^G = 0$, *in all degrees.* $H_S^*(X; \mathcal{O})^G = 0$, *unless $S$ is open.*

(b) *If $h$ is large enough (see* 2.10; $h > 0$ *suffices if $G$ is a torus*), $H_S^*(X; \mathcal{L}^h \otimes \mathcal{K})^G = 0$, *in all degrees.*

(2.7) *Remarks.* (a) $S$ can only be open if $X^{ss} = \emptyset$. If so, $H^*(S; \mathcal{O})^G = H^*(X; \mathcal{O})^G$.

(b) Part (b) can fail if $h$ is small and $G$ is not a torus; see examples (4.2) and (4.3).

*Proof.* $\mathbb{T}$ has no negative weights on $\Lambda^\bullet(\mathfrak{u})^t$, so that vanishing of (2.4), and then of (2.1), is guaranteed whenever the total $\mathbb{T}$-action on (2.5) has positive weights only. The weights of $(T_Z Y)^t$ and of $T_S X$ are positive, while $\det(T_S X)$ is $\mathbb{T}$-positive, unless $S$ is open. When $\mathcal{V}$ is $\mathcal{L}$ or $\mathcal{O}$, its $\mathbb{T}$-weights over $Z^\circ$ are also nonnegative, proving (a).

When $\mathcal{V} = \mathcal{L}^h \otimes \mathcal{K}$, we can factor

$$(2.8) \qquad \mathcal{K} = \mathcal{K}_Z \otimes \det^{-1}(T_Z Y) \otimes \det^{-1}(T_Y X) ,$$

where $\mathcal{K}_Z$ is the canonical bundle of $Z$. Note that, over $Y$, $T_S X = T_Y X / \bar{\mathfrak{u}}$, so (2.5) becomes

$$(2.9) \qquad \mathcal{L}^h \otimes \mathcal{K}_Z \otimes \mathrm{Sym}^\bullet (T_Z Y)^t \otimes \det^{-1}(T_Z Y) \otimes \mathrm{Sym}^\bullet(T_S X) \otimes \det^{-1}(\bar{\mathfrak{u}}) .$$

In the torus case, $\mathfrak{u} = 0$, so no $\mathbb{T}$-weights are negative, while the one on $\mathcal{L}^h$ is positive. In general, the problem factor is $\det \mathfrak{u}$. Subject to the restriction on $h$ spelled out below, we get part (b). □

(2.10) *Remark.* "Large enough" means that the $\mathbb{T}$-weight on $\mathcal{L}^h$ over $Z$ exceeds the one on $\det \bar{\mathfrak{u}} \otimes \det(T_Z Y)$. A sufficient condition is $h \cdot \|\beta\|^2 > \langle \beta | 2\rho \rangle$, the $\mathbb{T}$-weight on $\det \bar{\mathfrak{u}}$. A finer test, taking advantage of $T_Z Y$, replaces $\langle \beta | 2\rho \rangle$ by the $\mathbb{T}$-weight on $\det \bar{\mathfrak{s}}$, where $\mathfrak{s}$ is the part of $\mathfrak{u}$ fixing all points of $Z$. Note, by considering $\mathbb{T}$-weights, that a $\xi \in \mathfrak{u}$ vanishing at a point $z \in Z^\circ$ must vanish everywhere on $\varphi^{-1}(z)$; it follows in particular that any $h > 0$ will do, if the $\mathfrak{u}$-action on $Y$ is generically free.

(2.11) PROPOSITION. (a) $H^*(X; \mathcal{L}^h)^G = H^*(X^{\mathrm{ss}}; \mathcal{L}^h)^G$ for $h \geq 0$ (but we need $h > 0$ if $X^{\mathrm{ss}} = \emptyset$).
(b) $H^*(X; \mathcal{L}^{-h})^G = H^*_c(X^{\mathrm{ss}}; \mathcal{L}^{-h})^G$ when $h$ is large enough (cf. 2.10).

*Proof.* Part (a) follows from (2.6.a) and the ensuing collapse of the Cousin spectral sequence (A.6). The negative case follows from (2.6.b) and Serre duality (A.9), applied to $X$ and $X^{\mathrm{ss}}$. □

(2.12) *Remarks.* (i) As the proof shows, (a) applies to all vector bundles adapted to the stratification (strictly adapted, if $X^{\mathrm{ss}} = \emptyset$). This strengthens [TZ, Thm. 4.2], for holomorphic vector bundles.
(ii) From Kodaira's theorem, $H^*(X; \mathcal{L}^{-h}) = 0$, except in top degree. Subject to a condition (3.6), we shall see the same about the cohomology with proper supports $H^*_c(X^{\mathrm{ss}}; \mathcal{L}^{-h})^G$ (Serre duality in Prop. 6.2). Even then, the nonzero dimensions may disagree for small $h$; see (4.3). The difference can be computed from the extra terms in the Cousin sequence.

(2.13) *Remark.* (2.6) and (2.11) also hold when $X$ has *rational singularities*.[4] For an equivariant resolution $\pi : X' \to X$, we have $H_S^*(X; \mathcal{L}) = H_{\pi^{-1}(S)}^*(X'; \pi^*\mathcal{L})$, and similarly after twisting by Grauert's canonical sheaves [GR]. By (1.2), $\pi^{-1}(S)$ is a union of strata in $X'$, and $\pi^*\mathcal{L}$ has positive $\mathbb{T}$-weights, if $S$ is unstable; so the invariant part of the second space vanishes. For general singularities, (a) can fail for small $h$; see (4.6). For large $h$, we only have $H^0$, and then (2.11.a) holds whenever $X$ is normal.

## 3. Relation to quotient spaces

*Refresher on* GIT *quotients* [MFK, Ch. 1]. The GIT quotient $X//G$ is the projective variety $\operatorname{Proj} \bigoplus_{N \geq 0} \Gamma(X; \mathcal{L}^N)^G$. It is the scheme-theoretic quotient $X^{ss}/G$ of $X^{ss}$ by $G$, and the structural morphism $q : X^{ss} \to X//G$ is affine. The inclusion of $G$-invariants within $\operatorname{Spec} \bigoplus_{N \geq 0} \Gamma(X; \mathcal{L}^N)$ extends $q$ to a morphism between the affine cones over $X$ and $X//G$, mapping unstable points to the origin. *Stable* points have closed orbits and finite stabilizers. All semistable points are stable precisely when the $G$-action on $X^{ss}$ is proper; this happens if and only if all stabilizers are finite. If so, $X_G^{ss}$ is a compact, Kähler, Deligne-Mumford stack (A.8).

The *invariant direct image* $q_*^G \mathcal{F}$ of a (quasi) coherent, $G$-equivariant sheaf $\mathcal{F}$ over $X^{ss}$ (A.1) is the (quasi) coherent sheaf on $X//G$, whose sections over $U$ are the $G$-invariants in $\mathcal{F}(q^{-1}U)$. The functor $q_*^G$ is best viewed as the direct image along the morphism $q^G$, induced by $q$, from the stack $X_G^{ss}$ to $X//G$. As $q$ is affine and $G$ is reductive, $q_*^G$ is exact. The lift to $X^{ss}$ of a sheaf on $X//G$ has a natural $G$-structure; therefore, $(q^*, q_*^G)$ forms an adjoint pair, relating equivariant sheaves on $X^{ss}$ to sheaves on $X//G$. Further, $q_*^G \circ q^* = \operatorname{Id}$; in particular, $q_*^G \mathcal{O} = \mathcal{O}$.

An equivariant sheaf $\mathcal{F}$ over $X^{ss}$ *descends to* $X//G$ if it is $G$-isomorphic to a lift from downstairs. This happens if and only if adjunction $q^* q_*^G \mathcal{F} \to \mathcal{F}$ is an isomorphism; if so, $\mathcal{F} = q^* q_*^G \mathcal{F}$, and we shall abusively denote $q_*^G \mathcal{F}$ by $\mathcal{F}$ as well. Vector bundles descend precisely when the isotropies of closed orbits in $X^{ss}$ act trivially on the fibers (*Kempf's descent lemma*). Some power of $\mathcal{L}$ always descends, because the infinitesimal isotropies in $X^{ss}$ act trivially. Vector bundles with this last property are said to *descend fractionally*.

---

[4]Recall that a normal variety has rational singularities if and only if the higher direct images of $\mathcal{O}$, from any desingularization, vanish. Equivalently, it is Cohen-Macaulay and *Grauert's canonical sheaf* [GR] of completely regular top differentials (the push-down of the dualizing sheaf, from any desingularization) agrees with the *Grothendieck dualizing sheaf* of all top differentials.

(3.1) *Remark.* Affineness of $q$ and $q_*^G \mathcal{O} = \mathcal{O}$ are Seshadri's [S1] defining conditions for a *good quotient* under a reductive group action; the other properties of $q_*^G$ follow (in characteristic 0).

*Quantization commutes with reduction.* If $\mathcal{L}$ descends, $\mathbb{R}q_*^G \mathcal{L} = q_*^G \mathcal{L} = \mathcal{L}$ and (2.11.a) gives the following.

(3.2.a) THEOREM. *If $\mathcal{L}$ descends, $H^*(X;\mathcal{L})^G = H^*(X//G;\mathcal{L})$. If $X//G \neq \emptyset$, $H^*(X;\mathcal{O})^G = H^*(X//G;\mathcal{O})$.*

(3.3) *Remarks.* (i) When $\mathcal{L}$ does not descend, equality holds with $q_*^G \mathcal{L}$ downstairs. This need not be a line bundle, but, if nonzero, it is a reflexive sheaf of rank 1.

(ii) The theorem holds for any vector bundle adapted to the stratification (strictly adapted, if $X^{ss} = \emptyset$). Further, as in (2.13), rational singularities are permissible on $X$.

Negative powers of $\mathcal{L}$ raise a delicate question. According to Boutot [B], $X//G$ has rational singularities, so we can use Serre duality, once the invariant direct image of the canonical sheaf $\mathcal{K}$ of $X$ is known. We shall review that in a moment; most relevant is the *canonical sheaf $\mathcal{K}_G$ of the stack $X_G$*, the twist of $\mathcal{K}$ by the 1-dimensional representation $\det \mathfrak{g}$. (This plays the role of dualizing sheaf of $BG$; it is a sign representation, trivial when $G$ is connected.) Call $\omega$ the dualizing sheaf of $X//G$ and $\delta := \dim X - \dim X//G$.

(3.2.b) THEOREM. *If $\mathcal{L}^h$ descends and $h > 0$ is large enough (2.10), an isomorphism $q_*^G \mathcal{K}_G \cong \omega$ determines another one:*

$$\left[ H^*(X;\mathcal{L}^{-h}) \otimes \det \mathfrak{g} \right]^G \cong H^{*-\delta}(X//G;\mathcal{L}^{-h}) .$$

*These vanish if $* \neq \dim X$. If stable points exist and the $\mathfrak{g}$-action on $X$ is free in codimension 1, any $h > 0$ will do, and $q_*^G \mathcal{K}_G = \omega$, canonically.*

*Proof.* When $\mathcal{L}^h$ descends,

$$q_*^G(\mathcal{L}^h \otimes \mathcal{K}_G) = \mathcal{L}^h \otimes q_*^G \mathcal{K}_G \cong \mathcal{L}^h \otimes \omega .$$

From (2.11.b) and Serre duality,

(3.4)
$$\left[ H^{\text{top}}\left(X;\mathcal{L}^{-h}\right) \otimes \det \mathfrak{g} \right]^G = H_G^0\left(X;\mathcal{K}_G \otimes \mathcal{L}^h\right)^t = H_G^0\left(X^{ss};\mathcal{K}_G \otimes \mathcal{L}^h\right)^t$$
$$\cong H^0\left(X//G;\mathcal{L}^h \otimes \omega\right)^t = H^{\text{top}}\left(X//G;\mathcal{L}^{-h}\right) .$$

Cohomology vanishing is Kodaira's theorem. The last part of (3.2.b) follows from the criterion in (2.10) (strata of codimension $\geq 2$ do not affect $H^0$) and from Knop's theorem (3.6) below. □

(3.5) *Remarks.* (i) Call $\mathbb{R}q_!^G$ the derived invariant direct image with proper supports along $q$, shifted down by $\delta$ (A.9). Under relative duality, an isomorphism $q_*^G \mathcal{K}_G \cong \omega$ correspond to a quasi-isomorphism $\mathbb{R}q_!^G \mathcal{O} \cong \mathcal{O}$, and the theorem also follows by application of $\mathbb{R}q_!^G$ in (2.11.b).

(ii) Irrespective of $q_*^G \mathcal{K}_G$, (3.4) holds, for large enough $h$, with $\mathbb{R}q_!^G \mathcal{L}^{-h}$ downstairs. Freedom of the action in codimension 1 ensures that $\mathbb{R}q_!^G \mathcal{L}^{-h} = q_*^G \mathcal{L}^{-h}$, even when $\mathcal{L}^h$ does not descend, and that the second isomorphism holds for any $h > 0$.

*Knop's results on dualizing sheaves.* Keep $X$ smooth, although Knop's results [Kn] hold more generally. Consider the condition:

(3.6) Stable points exist, and the $G$-action on $X^{ss}$ has finite stabilizers in codimension 1.

(3.7) (Cf. [Kn, Kor. 2].) *If* (3.6) *holds, there exists a natural isomorphism* $\psi : q_*^G \mathcal{K}_G \to \omega$.

A description of $\psi$ different from Knop's will be needed in Section 6, so we shall reprove this in Section 5. Note, meanwhile, a more general result. The kernel $\mathfrak{k}$ of the infinitesimal action $\mathfrak{g} \otimes \mathcal{O} \to TX$ is a vector bundle in codimension 1; call $\lambda$ the line bundle extension of $\det^{-1}\mathfrak{k}$ to all of $X^{ss}$.

(3.8) [Kn, Kor. 1]. *If the generic fiber of $q$ contains a dense orbit, then, for an effective $G$-divisor $D$ supported by the points where the stabilizer dimension jumps, $\omega$ is naturally isomorphic to $q_*^G(\lambda \otimes \mathcal{K}_G(D))$.*

Note that $\lambda$ corrects for generic positive-dimensional stabilizers, but I do not know a good interpretation for $D$. Note also that, when the generic orbits in $X^{ss}$ are closed, a jump in stabilizer dimension in codimension 1 entails the appearance of a unipotent radical in the isotropy (use a slice argument).

(3.9) COROLLARY. *If the generic orbits in $X^{ss}$ are closed, then, for large $h$, $H^*\left(X//G; q_*^G(\mathcal{L}^{-h})\right)$ is the $\det^{-1}\mathfrak{g}$-typical component of*

$$H^{*+\delta}\left(X; \lambda^{-1} \otimes \mathcal{L}^{-h}(-D)\right) .$$

*Generic points in $D$ have nonreductive isotropy.*

*Addendum*: *Shifted quotients.* Meinrenken and Sjamaar [MS, Thm. 2.14] prove a different statement for $\mathcal{L}^{-1}$, involving *shifted quotient* $(X \times F_{2\rho})//G$, linearized by $\mathcal{L}(2\rho)$ (cf. 1.3.iii). For convenience, we call it $X_{2\rho}//G$. This cannot be easily described in terms of stacks ($X_{2\rho}//G$ depends on $X$ and $G$, not only on the stack $X_G$), but has the advantage of removing the "large $h$" restriction in (3.2.b). Further, since all semistable points in $X \times F_{2\rho}$ have reductive isotropies, (3.6) is now replaced by the simpler condition that the shifted quotient should contain stable points. At the referee's suggestion, I shall sketch a proof of their result along the lines of Section 2. (Due to Kodaira vanishing, this does not strengthen the Riemann-Roch statement of [MS], except when $X_{2\rho}//G$ is singular.) We also note amusing half-way versions, involving $\mathcal{K}^{1/2}$ and the $\rho$-shifted quotient.

(3.10) PROPOSITION.   (a) *If $X_{2\rho}//G$ contains stable points,*

$$H^{\mathrm{top}}\left(X; \mathcal{L}^{-1}\right)^G = H^{\mathrm{top}}\left(X_{2\rho}//G; q_*^G(\mathcal{L}(2\rho)^{-1})\right),$$

*and the other cohomologies vanish.*
(b) *Always,*

$$H^*\left(X; \mathcal{K}_X^{1/2} \otimes \mathcal{L}\right)^G = H^*\left(X_\rho//G; q_*^G\left(\mathcal{K}_{X\times F}^{1/2} \otimes \mathcal{L}(\rho)\right)\right).$$

(c) *If $X_\rho//G$ contains stable points,*

$$H^*\left(X; \mathcal{K}_X^{1/2} \otimes \mathcal{L}^{-1}\right)^G = H^{*-\delta}\left(X_\rho//G; q_*^G\left(\mathcal{K}_{X\times F}^{1/2} \otimes \mathcal{L}(\rho)^{-1}\right)\right).$$

(3.11) *Remarks.* (i) If $G$ is not connected, we must twist $\mathcal{K}$ by $\det(\mathfrak{g})$.
(ii) It will emerge from the proof that the conditions on stable points in (a) and (c) can be much weakened, if one is interested in $(X \times F)^{\mathrm{ss}}$ and not in the quotient.

*Proof (sketch).* To argue as in (3.2.b), we must show that $H_S^*(X \times F; \mathcal{K} \otimes \mathcal{L}(2\rho))^G = 0$, for an $\mathcal{L}(2\rho)$-unstable stratum $S$ on $X \times F$, with associated $\beta$ and $Z = Z_X \times Z_F$. If $W_L$ and $W_G$ denote the Weyl groups of $L$ and $G$, the $H$-fixed point set of $Z_F$ is $W_L \cdot w$, for some $w \in W_G$. Existence of $L$-semistable points for $\mathcal{L}(2\rho - \beta)$ on $Z$ requires $\alpha := \beta - 2w\rho$ to belong to the $H$-moment map image of $Z_X$, and this in turn forces the nilradical of the generic isotropy on $Z_X$ to lie in the $\alpha$-negative parabolic $\mathfrak{q} \subseteq \mathfrak{g}$. The $\mathbb{T}$-negative part $\mathfrak{s}$ of the generic isotropy of $Z^\circ$ (cf. 2.10) lies then in $\mathfrak{q} \cap w\mathfrak{n}$. This is $v\mathfrak{n} \cap w\mathfrak{n}$, for a certain $v \in W_G$ for which $v^{-1}\alpha$ is dominant; so $\langle \beta | \det \mathfrak{s} \rangle$ is underestimated by $\langle \beta | \det(v\mathfrak{n} \cap w\mathfrak{n}) \rangle = -\langle \beta | v\rho + w\rho \rangle$, and we have

(3.12)
$$\begin{aligned}
\beta^2 + \langle \beta | \det \mathfrak{s} \rangle &\geq \langle \beta | \beta - v\rho - w\rho \rangle \\
&= (\beta - v\rho - w\rho)^2 + \langle v\rho + w\rho | \beta - v\rho - w\rho \rangle \\
&= (\beta - v\rho - w\rho)^2 + \langle v\rho + w\rho | \alpha \rangle + \langle v\rho + w\rho | w\rho - v\rho \rangle \\
&= (\beta - v\rho - w\rho)^2 + \langle v\rho + w\rho | \alpha \rangle.
\end{aligned}$$

The second term is semipositive ($v^{-1}\alpha$ is dominant), so that (3.12) is positive, unless $\beta = v\rho + w\rho$, $\alpha = v\rho - w\rho$. In the latter case, dominance of $v^{-1}\alpha$ implies that $\alpha = \det \bar{\mathfrak{q}}$, $\beta = 2\rho_M$, where $M$ is the Levi component of $\mathfrak{q}$. The longest element $w' \in W_M \subset W_G$ fixes $\alpha$ and changes the sign of $\beta$; thus it acts on $Z^\circ$, swapping positive and negative normal directions. If $Y^\circ$ is the $N$-orbit of $Z^\circ$, the action of $w'$ shows that $S = GY^\circ$ is open in $X$; in this case, $X_\rho//G$ is empty (unless $\beta = 0$, but then $S$ was semistable). Else, there must be extra directions in $T_Z Y$, not yet counted in the estimate (3.12) of $\det^{-1}(T_Z Y)$; and they force vanishing of $\mathbb{T}$-invariants in (2.9).

The proof of (3.10.b) is similar, but involves $\rho/2$ and $\det^{1/2}\mathfrak{s}$. □

## 4. Counterexamples

Occasional failure of the naive quantization statement for negative bundles was already noted by M. Vergne [MS, §2.14]. The examples below justify the restrictions in (3.2.b) and other theorems.

(4.1) The statement can even fail for $G = \mathbb{C}^\times$, in the *absence* of stable points. Let $X = \mathbb{P}^1$, $\mathcal{L} = \mathcal{O}(1)$, with the obvious $\mathbb{C}^\times$-action on $\mathbb{P}^1$ and its lift to $\mathcal{L}$ which fixes the fiber over $0$. $X^{ss} = \mathbb{P}^{-1} - \{\infty\}$, $X//G$ is a point, but $q_*^G \mathcal{K} = 0$. The invariant cohomology of negative powers of $\mathcal{L}$ vanishes, yet the space of sections on the quotient point is always a line. The correction divisor $D$ in (3.7) is the origin, with multiplicity $2$.

(4.2) For an example where $H_S^*(X; \mathcal{K} \otimes \mathcal{L})^G \neq 0$, let $G$ be simple, $B \subset G$ a Borel subgroup, $X = G/B \times G/B$, with $G$ acting diagonally. With $\mathcal{L} = \mathcal{O}(\rho) \boxtimes \mathcal{O}(\rho)$, $\mathcal{K} \otimes \mathcal{L} = \mathcal{O}(-\rho) \boxtimes \mathcal{O}(-\rho)$. The strata are the $G$-orbits, labeled by Weyl group elements; they correspond to the $B$-orbits on $G/B$. Also, $H_S^{\mathrm{codim}(S)}(X; \mathcal{K} \otimes \mathcal{L})^G = \mathbb{C}$ for each orbit $S$, while the other cohomologies vanish. The Cousin sequence collapses at $E_2$, not $E_1$.

(4.3) Another counterexample to the statement for small negative powers, where stable points exist, arises from the multiplication action of $SL_2$ on $X = \mathbb{P}^3$, the projective space of $2 \times 2$ matrices, and $\mathcal{L} = \mathcal{O}(1)$. The unstable locus is the quadric surface of singular matrices, while $X^{ss}$ is a single orbit isomorphic to $\mathbb{PSL}_2$. The stabilizer $\mathbb{Z}/2$ acts nontrivially on $\mathcal{L}$. Serre duality identifies $H_c^3(\mathbb{PSL}_2; \mathcal{O}(-h))$, the only nonzero group, with the dual of $H^0(\mathbb{PSL}_2; \mathcal{O}(h-4))$. The invariants vanish for odd $h$, and equal $\mathbb{C}$ when $h > 0$ is even. This matches $H^*(\mathbb{P}^3; \mathcal{O}(-h))^{SL_2}$ for all $h > 0$, except $h = 2$.

(4.4) Without the codimension condition in (3.6), isomorphism in (3.2.b) can fail for all negative powers of $\mathcal{L}$. Consider the multiplication action of $SL_2$ on the space of $2 \times 2$ matrices, completed to $\mathbb{P}^4$ by addition of a hyperplane at infinity. The quotient is $\mathbb{P}^1$. The unstable set is the quadric surface of singular matrices at infinity. Since $q^*\mathcal{O}(1) = \mathcal{O}(2)$, $q_*^G \mathcal{O}(2n) = \mathcal{O}(n)$ for all $n$.

One sees on global sections that $q_*^G \mathcal{O}(1) = \mathcal{O}$. Thus, $q_*^G \mathcal{O}(2n+1) = \mathcal{O}(n)$, and $q_*^G \mathcal{K} = q_*^G \mathcal{O}(-5) = \mathcal{O}(-3) \neq \omega$. Quantization fails to commute with reduction for nearly all negative bundles, even on Riemann-Roch numbers: $h^1(\mathbb{P}^1; \mathcal{O}(-h)) = h - 1$, and $h^0 = 0$; but the invariants in $H^*(\mathbb{P}^4; \mathcal{O}(-2h))$ are in degree 2 and have dimension $h - 2$ (if $h \geq 2$). A twist by $\mathcal{O}(-2)$ on $X = \mathbb{P}^4$ corrects the statement; $D$, in (3.7), is the cone of singular matrices.

(4.5) An interesting example (without stable points) is the adjoint representation of $SL_2$, with quotient $\mathbb{A}^1$. Completing $\mathfrak{sl}_2$ to $\mathbb{P}^3$, the GIT quotient becomes $\mathbb{P}^1$, the unstable locus being the trace of the nilpotent cone at infinity. Again, $\mathcal{O}(-2) = q^*\mathcal{O}(-1)$, but this time $\mathcal{K}$ pushes down correctly: $q_*^G \mathcal{O}(-4) = \mathcal{O}(-2)$. Consider now $\tilde{X}$, the blow-up of the origin, with exceptional divisor $E$. Linearizing by $\mathcal{O}(2)(-E)$, the quotient is still $\mathbb{P}^1$ (the unstable locus is the proper transform of the nilpotent cone). However, the push-down of the canonical sheaf is now $\mathcal{O}(-1)$. The natural morphism $q_*^G \mathcal{K} \to \tilde{q}_*^G \tilde{\mathcal{K}}$, which normally gives rise to $\psi$ in (3.6), has a 1-dimensional cokernel at 0. The problem stems not from a jump in the stabilizer dimension, but from the non-trivial dualizing line of the (dihedral) stabilizer $\mathbb{C}^\times \tilde{\times} \mathbb{Z}/2$ along $E$. (In (3.7), $\lambda$ is nontrivial.)

(4.6) Finally, (3.2.a) can fail if $X$ has irrational singularities. Choose a smooth proper curve $\Sigma$ of genus $\geq 2$ and a positive line bundle $\mathcal{F}$ with $H^1(\Sigma; \mathcal{F}) \neq 0$. Let $X$ be the cone over $\Sigma$, obtained by adding a point to the total space of $\mathcal{F}$, $\mathcal{L}$ the tautological line bundle on $X$ restricting to $\mathcal{F}$ over $\Sigma$ and to $\mathcal{O}(1)$ on every generator. Lift the $\mathbb{C}^\times$-action on $X$ to $\mathcal{L}$, making it trivial over $\Sigma$. The unstable locus is the vertex, and the GIT quotient is the original curve, over which $\mathcal{L} = \mathcal{F}$ has nontrivial first cohomology. Yet $H^1(X; \mathcal{L}) = 0$ (say, by a Mayer-Vietoris calculation), while $H_G^2$ with supports at the vertex equals $H^1(\Sigma; \mathcal{F})$.

## 5. Relative version, rigidity and rational singularities

This section presents the relative version (5.2) of the quantization theorem, whose special case (5.4) has the attractive applications (5.5), (5.6), and Boutot's theorem (5.7). The more technical wall-crossing lemma (5.8) will be used for a vanishing theorem in Section 6.

*Relative version of the theorem.* Recall, in the discussion (1.2) of a projective morphism $\pi : X' \to X$, that $(X')^{ss} \subseteq \pi^{-1} X^{ss}$; thus, $\pi$ induces a map $p : X'//G \to X//G$ on quotients. The former is $\text{Proj} \bigoplus_{n \geq 0} q_*^G \pi_* \mathcal{M}^n$ over the latter, whereas $\pi^{-1} X^{ss}$ is $\text{Proj} \bigoplus_{n \geq 0} q_* \pi_* \mathcal{M}^n$. The difference $\pi^{-1}(X^{ss}) - (X')^{ss}$ is the base locus of invariant relative sections of powers of $\mathcal{M}$, and $q' : (X')^{ss} \to X'//G$ arises from the obvious inclusion of graded algebras over $X//G$.

(5.1) LEMMA. *The $\mathcal{L}_\varepsilon$-stratification ($\varepsilon \to 0$) on $\pi^{-1}(X^{\mathrm{ss}})$ depends on $\mathcal{M}$ and $X^{\mathrm{ss}}$, but not on $\mathcal{L}$.*

*Proof.* The description just given for $(X')^{\mathrm{ss}}$ does not involve $\mathcal{L}$. Further, any $\beta$ and $\mathbb{T}$ in $\pi^{-1}(X^{\mathrm{ss}})$ depend on the $\mathcal{M}$-weights alone (the $\mathcal{L}$-weights are null), and so we can argue inductively for the other strata. □

This allows us to stratify $X'$ when $X$ is *affine*, taking $\mathcal{L} = \mathcal{O}$, with affine quotient $X/\!/G = X/G$. The construction is local on $X/G$ ($G$-local near *closed* orbits in $X = X^{\mathrm{ss}}$). The stratification can be described word-for-word as in Section 1, except the $Z^\circ/L$ are now projective over $X/G$. In the rest of this section, either $X$ is projective and $\mathcal{L}$ is ample, or else $X$ is affine and $\mathcal{L} = \mathcal{O}$.

(5.2) RELATIVE QUANTIZATION. *Let $X'$ have rational singularities. If $h \geq 0$ and the $\mathcal{M}^h$ descend to $X'/\!/G$, $\mathbb{R}p_*\mathcal{M}^h = q^G_*\mathbb{R}\pi_*\mathcal{M}^h$ on $X/\!/G$ ($h > 0$ is needed if any $X'$-component has empty quotient).*

(5.3) *Remarks.* (i) If the $G$-action on $X$ is trivial, this becomes $\mathbb{R}p_*\mathcal{M}^h = (\mathbb{R}\pi_*\mathcal{M}^h)^G$, which is a quantization theorem for a family parametrized by $X$. Taking $X$ to be a point simply recovers (3.2.a) for $X'$ and $\mathcal{M}$.

(ii) More generally, if less simply,
$$\mathbb{R}p_* \circ q'^G_*(\mathcal{M}^h \otimes \pi^*\mathcal{V}) = q^G_*(\mathcal{V} \otimes \mathbb{R}\pi_*\mathcal{M}^h),$$
for any $h \geq 0$ and $G$-vector bundle $\mathcal{V}$ on $X^{\mathrm{ss}}$ which descends fractionally (with the usual caveat for $h = 0$).

(5.4) COROLLARY (assumptions as in 5.2). *If $\mathbb{R}\pi_*\mathcal{O} = \mathcal{O}$, then $\mathbb{R}p_*\mathcal{O} = \mathcal{O}$; and, for any vector bundle $\mathcal{W}$ on $X/\!/G$, $H^*(X'/\!/G; p^*\mathcal{W})$ and $H^*(X/\!/G; \mathcal{W})$ are naturally isomorphic.*

Below, we shall deduce Boutot's theorem (5.7) from (5.4). The two are, in fact, equivalent. (However, (5.2) is strictly stronger; cf. the proof of 5.8.) Indeed, the assumptions in (5.4) ensure, possibly after we replace $\mathcal{M}$ by a large power, that $\mathrm{Spec}\bigoplus_{n \geq 0} \Gamma(X'; \mathcal{M}^n)$ has rational singularities. So, also, does the quotient $\mathrm{Spec}\bigoplus_{n \geq 0} \Gamma(X'; \mathcal{M}^n)^G$; and then, $\mathcal{O}$ has no higher cohomology over $X'/\!/G = \mathrm{Proj}\bigoplus_{n \geq 0} \Gamma(X'; \mathcal{M}^n)^G$. This argument is due to Broer and Sjamaar in the absolute case (when $X$ is a point) [Sj1, 2.23]. Curiously, I do not find the relative case and its consequences (5.5), (5.6) in the literature.

*Proof of* (5.2). For affine $X$, this amounts to
$$\Gamma\left(X/\!/G; q^G_*\mathbb{R}^i\pi_*\mathcal{M}^h\right) = H^i\left(X'/\!/G; \mathcal{M}^h\right).$$

Lifting to $X'$, this means

$$H_G^i\left(\pi^{-1}X^{\mathrm{ss}};\mathcal{M}^h\right) = H_G^i\left((X')^{\mathrm{ss}};\mathcal{M}^h\right),$$

which follows from vanishing of the invariant cohomology of $\mathcal{M}^h$ with supports on $\pi^{-1}(X^{\mathrm{ss}}) - (X')^{\mathrm{ss}}$. The affine case suffices, as the statement is local on $X/\!/G$. □

*Application*: *Rigidity and rational singularities of quotients.* Important applications of (5.4) arise from the examples (1.3), *when $X$ itself has rational singularities.* The first two results allow us to replace badly singular quotients with nicer ones, without changing the cohomology of vector bundles. Perturbing the polarization of $X$ slightly gives, for a vector bundle $\mathcal{W}$ on $X/\!/G$:

(5.5) RIGIDITY THEOREM. *If the perturbed quotient $X_\varepsilon/\!/G$ is not empty,*

$$H^*\left(X_\varepsilon/\!/G; p^*\mathcal{W}\right) = H^*\left(X/\!/G; \mathcal{W}\right).$$

When no perturbation leads to a nice quotient, shift desingularization (1.3.iii) can help. Call $X'_\lambda/\!/G$ the GIT quotient of $X \times F_\lambda$ under the diagonal action of $G$, linearized by $\mathcal{L} \boxtimes (\lambda)$ (for small fractional $\lambda$).

(5.6) SHIFTING THEOREM. *A nonempty $X'_\lambda/\!/G$, with small $\lambda$, can replace $X'/\!/G$ in (5.5).*

(5.7) BOUTOT'S THEOREM. *$X/\!/G$ has rational singularities, if $X$ does.*

*Proof.* If $X$ is singular, choose a resolution $X'$. Next, recall Kirwan's partial desingularization of $X'/\!/G$ by the GIT quotient $\tilde{X}/\!/G$ of a sequential blow-up $\tilde{X}$ of $X'$ along smooth $G$-subvarieties [K2]: at every stage, the center of blowing-up comprises the points whose stabilizers have maximal dimension. $\tilde{X}/\!/G$ has finite quotient singularities, which are rational (Burns [Bu], Viehweg [Vi]; see also Remark 5.11). If $\tilde{X}/\!/G \neq \emptyset$, the composite $p : \tilde{X}/\!/G \to X/\!/G$ is birational. From (5.4) and the assumption on $X$, $\mathbb{R}p_*\mathcal{O} = \mathcal{O}$.

This argument fails if $X'$ or $\tilde{X}$ lead to an empty quotient. By [K2, 3.11], this only happens if, at some stage of the process, the center of blowing-up dominates the quotient. In that case, we can replace $X$ by that (smooth) blowing-up center, without changing the quotient, and proceed as before. □

*Application*: *canonical wall-crossing.* When $X = X'$, two quotients $X_+^{\mathrm{ss}}/G$, $X_-^{\mathrm{ss}}/G$ (for different signs of $\varepsilon$) can be compared with $X^{\mathrm{ss}}/G$, and we may consider vector bundles which descend to the perturbed quotients, but not to $X^{\mathrm{ss}}/G$. Showing equality of the plus and minus cohomologies (invariance under wall-crossing) by the argument in Section 2 requires more information about the strata. (Typically, we need absolute bounds on the fiberwise $\mathbb{T}$-weights;

the bundle should not be "too far" from descending.)[5] An important exception concerns $\mathcal{M} = \mathcal{K}$. (C. Woodward informs me that the result was known for Riemann-Roch numbers of smooth symplectic quotients.)

(5.8) LEMMA. *If $X^{ss}$ is smooth and verifies* (3.6), *$\mathcal{V}$ descends fractionally to $X^{ss}/G$ and $0 \leq \alpha < 1$, the cohomologies of the invariant direct images of $\mathcal{V} \otimes \mathcal{K}^\alpha$ on all three quotients agree.*

*Proof.* Equality for $X_+//G$ and $X//G$ is an instance of (5.3.ii), without any upper bound on $\alpha$. On the minus side, $\mathcal{V} \otimes \mathcal{K}^\alpha = (\mathcal{V} \otimes \mathcal{K}^{\alpha-1}) \otimes \mathcal{K}$, and $\mathcal{K}^{\alpha-1}$ is positive relative to $p_- : X_-//G \to X//G$. A classical argument of Grauert and Riemenschneider [GR, Satz 2.3] deduces the vanishing of all higher direct images of $q^G_{-*}(\mathcal{V} \otimes \mathcal{K}^\alpha)$ under $p_-$ from a global vanishing result on $X_-//G$; that vanishing (which slightly extends [GR, Satz 2.2]) is the subject of Theorem 6.2 below. Agreement of the direct images requires elimination of cohomologies supported on codimension 1 strata. This follows from the last criterion in (2.10) and the assumed freedom of the $\mathfrak{g}$-action (thus of the $\mathfrak{u}$-action) in codimension 1. □

*Appendix*: *Construction of $\psi$ and proof of* (3.7).

*Stable case.* The result is clear when $G$ acts freely on $X^{ss}$; given that, étale-locally on $X//G$, Luna's slice theorem (A.7) reduces the case of a proper action to a finite group quotient. Assume then that $G$ and $q$ are finite; $q_!\mathcal{O} = q_*\mathcal{O}$ splits, as an $\mathcal{O}_{X//G}$-module, into $G$-isotypical components, with $\mathcal{O}_{X//G}$ as the identity summand. Grothendieck duality along $q$ gives natural isomorphisms (note that $\mathcal{K}_G$ is now $\mathcal{K}$):

$$(5.9) \quad \omega = \mathbb{R}\mathcal{H}om_{X//G}(\mathcal{O}; \omega) = \mathbb{R}\mathcal{H}om_{X//G}(q_!\mathcal{O}; \omega)^G$$
$$= \mathbb{R}q_* \circ \mathbb{R}\mathcal{H}om_{X//G}\left(\mathcal{O}; q^!\omega\right)^G = q^G_*\mathcal{K}_G.$$

(5.10) *Remark.* In stack language, we would say that the morphism $q^G : X^{ss}_G \to X//G$ is proper; thus, $q^G_* = q^G_!$ has a right adjoint $q^{G!}$, for which $q^{G!}\omega = \mathcal{K}_G$. Since $\mathbb{R}q^G_!\mathcal{O} = \mathcal{O}$, duality for $q^G$ gives

$$\omega = \mathbb{R}\mathcal{H}om_{X//G}(\mathcal{O}; \omega) = \mathbb{R}q^G_* \circ \mathbb{R}\mathcal{H}om_{X^{ss}_G}(\mathcal{O}; \mathcal{K}_G) = q^G_*\mathcal{K}_G.$$

---

[5]This is exactly right when the strata on the two sides of the perturbation are in one-to-one correspondence, as in a *truly faithful* wall-crossing, studied in [DH, §4] and [Th]. Each $\mathbb{T}$ is then the stabilizer of a $Z$, and the bounds on the fiberwise $\mathbb{T}$-weights over $Z$ come from the determinants of the positive and negative parts of the normal bundles; cf. (2.6).

(5.11) *Remark.* This shows that $X//G$ is Cohen-Macaulay, since the dualizing complex $q_*^G \mathcal{K}_G$ lives in a single degree. Now, for finite group quotients, it is clear that top differentials are completely regular (locally $L^2$ [GR, §2]) downstairs, if their liftings are so; thus, Grauert's and Grothedieck's canonical sheaves agree on the quotient, if they did so upstairs. Rational singularities are thus preserved by finite quotients.

*General case.* Let, for simplicity, $X^{\text{ss}}$ be smooth. In the proof of (5.7),

$$q_*^G \mathcal{K}_G = q_*^G \circ \pi_* \tilde{\mathcal{K}}_G .$$

With $\tilde{q} : \tilde{X}^{\text{ss}}//G$ denoting the structural map and $i : \tilde{X}^{\text{ss}} \subset \pi^{-1}(X^{\text{ss}})$ the inclusion, we have canonical isomorphisms

(5.12) $$\omega = p_* \tilde{\omega} = p_* \circ \tilde{q}_*^G \tilde{\mathcal{K}}_G = q_*^G \circ \pi_* \circ i_* \tilde{\mathcal{K}}_G .$$

The first holds because $X//G$ and $\tilde{X}//G$ have rational singularities; $\tilde{q}_*^G \tilde{\mathcal{K}}_G = \tilde{\omega}$ has just been proven, while $q \circ \pi \circ i = p \circ \tilde{q}$ gives the last equality. Also, $\psi : q_G^* \mathcal{K}_G \to \omega$ is defined from $\mathcal{K}_G = \pi_* \tilde{\mathcal{K}}_G \subset \pi_* \circ i_* \tilde{\mathcal{K}}_G$.

We check that $\psi$ is an isomorphism stepwise in Kirwan's construction. Luna's slice theorem and the result for finite groups reduce us to the case of a linear representation $X = V$ of a connected reductive group $R$, with $V^R = \{0\}$, and its blow-up $\tilde{X} = \tilde{V}$ about the origin. Let $\kappa$ be the canonical bundle of $\mathbb{P}(V)$. The desired isomorphism can be checked on global sections over $V$, and is guaranteed by the vanishing, for $n > 0$, of $R$-invariants in the cohomologies $H_S^1(\mathbb{P}(V); \kappa(n))$ supported on strata $S$ of codimension 1. This follows from the last criterion in (2.10); the $\mathfrak{u}$-action on $Y$ is generically free, since (3.6) requires this on the cone over $Y$ in $V$, and nilpotence of $\mathfrak{u}$ ensures its generic freedom on $Y \subset \mathbb{P}(V)$.

(5.13) *Remark.* Similarly, $q_*^G(\mathcal{V} \otimes \mathcal{K}) = p_* \circ \tilde{q}_*^G(\pi^* \mathcal{V} \otimes \tilde{K})$, if $\mathcal{V}$ descends fractionally.

## 6. Cohomology vanishing

Let $X^\circ$ be a smooth $G$-variety having a projective good quotient (3.1) (e.g. the earlier $X^{\text{ss}}$). Call an equivariant line bundle over $X^\circ$ *G-ample* if it descends fractionally to an ample line bundle on $X^\circ/G$. This is stronger than ampleness. (For instance, if $X^\circ$ is quasi-affine, ampleness is a vacuous condition.)

(6.1) LEMMA. *$X^\circ$ can be realized as the semistable part of a smooth, proper $\bar{X}$, carrying a linearization extending (some power of) any given G-ample line bundle on $X^\circ$.*

*Proof.*[6] We shall, in fact, produce a relative $G$-projectivization $\bar{q}: \bar{X} \to X^\circ/G$ of $q$, with relatively ample invariant divisor $E$ at infinity; and we shall see that $\bar{X}$ can be chosen smooth, if $X^\circ$ is so. Granting that, choose a large $\mathcal{L}$ on $X^\circ/G$; then, $\bar{q}^*\mathcal{L}(E)$ is ample on $\bar{X}$ and carries an obvious $G$-action. Thus linearized, the semistable part of $\bar{X}$ is precisely $X^\circ$, because all invariant sections of $\bar{q}^*\mathcal{L}^n$ over $X^\circ$ extend by zero to sections of $\bar{q}^*\mathcal{L}^n(nE)$ over $\bar{X}$.

To find $\bar{q}$, choose a coherent $G$-subsheaf $\mathcal{F} \subset q_*\mathcal{O}$, large enough to embed $X^\circ$ properly into the linear space $\Lambda := \operatorname{Spec}\operatorname{Sym}^\bullet(\mathcal{F})$ over $X^\circ/G$. The closure $\bar{X}$ of $X^\circ$ in the relative projective completion of $\Lambda$ and the relatively ample divisor $E$ at infinity are equivariant. If $X^\circ$ is smooth, $\bar{X}$ will become so, after sequentially blowing up $G$-subvarieties in $\bar{X} - X^\circ$; and the desired $E$ is obtained by subtracting, at each blow-up step, from the pull-back of the previous $E$, a small multiple of the exceptional divisor. $\square$

(6.2) PROPOSITION. *If* (3.6) *holds and $\mathcal{F}$ is $G$-ample, higher cohomology of $q_*^G(\mathcal{F} \otimes \mathcal{K}_G)$ over $X^\circ/G$ vanishes.* (*Equivalently, the invariant higher cohomology of $\mathcal{F} \otimes \mathcal{K}_G$ over $X^\circ$ vanishes.*)

*Proof.* The lifting $\tilde{\mathcal{F}}$ of $\mathcal{F}$ to $\tilde{X}^\circ$ (cf. the proof of 5.7) defines a quasi-positive line bundle over the smooth, compact, Kähler DM stack $\tilde{X}^\circ//G$. By Grauert-Riemenschneider [GR],
$$H^{>0}(\tilde{X}^\circ//G; \tilde{\mathcal{F}} \otimes \tilde{\mathcal{K}}) = 0$$
and
$$\mathbb{R}^{>0}p_* \circ \tilde{q}_*^G(\tilde{\mathcal{F}} \otimes \tilde{\mathcal{K}}_G) = 0 \ .$$
As in Remark (5.13), $p_* \circ \tilde{q}_*^G(\tilde{\mathcal{F}} \otimes \tilde{\mathcal{K}}_G) = q_*^G(\mathcal{F} \otimes \mathcal{K}_G)$. $\square$

(6.3) *Remark.* This can fail without (3.6): if $\mathcal{F} = \mathcal{O}(1)$ in (4.4),
$$q_*^G(\mathcal{K} \otimes \mathcal{F}) = \mathcal{O}(-2) \ ,$$
on $\mathbb{P}^1$, has $H^1 \neq 0$.

The results of Section 5 allow us to relax the assumptions. Let $X_+^\circ$ be the semistable locus for a perturbed linearization (or shift desingularization) of $X^\circ$; note, because of (6.1), that the refinement lemma (1.2) applies.

(6.4) PROPOSITION. *If $\mathcal{L}$ descends fractionally to $X^\circ/G$, $\mathcal{F} := \mathcal{L} \otimes \mathcal{K}_G^{-1}$ is $G$-ample on $X_+^\circ$, and the latter verifies* (3.6), *then* $H^{>0}\left(X^\circ/G; q_*^G\mathcal{L}\right) = 0$.

When $X$ is smooth and projective, and $\mathcal{L} > 0$ and $\mathcal{L} > \mathcal{K}$, Kodaira's theorem and (3.2.a) imply that $H^{>0}(X//G; \mathcal{L}) = 0$. Yet $\mathcal{L}$ need not dominate $\mathcal{K}$ on $X//G$, since the quotient with respect to $\mathcal{L}$ may differ from the quotient at $\mathcal{L}\mathcal{K}^{-1}$, where the latter is guaranteed to be ample. (The difference between quotients is often described in terms of $\mathcal{K}$-flips [Th], which create negative subvarieties for $\mathcal{L}\mathcal{K}^{-1}$.) The previous argument for vanishing breaks down when $X$ is not compact; and indeed, the conclusion can then fail (example 6.6). A useful criterion arises from the wall-crossing lemma. (Astute use of Nadel's vanishing theorem in (6.4) might well accomplish the same, but I have only checked special cases.)

---

[6]The author thanks Y. Hu for spotting an error in the original proof.

(6.5) PROPOSITION. $H^{>0}\left(X^{\mathrm{ss}}/G; q_*^G\mathcal{L}\right) = 0$, when the following hold:

(i) $\mathcal{L}$ and $\mathcal{L}\mathcal{K}^{-1}$ are ample on $X$;

(ii) all quotients $X_t^{\mathrm{ss}}/G$ (linearized by $\mathcal{L}\mathcal{K}^{-t}$) are projective;

(iii) the family $X_t^{\mathrm{ss}}$ is "upper semicontinuous" for $0 \le t \le 1$;

(iv) (3.6) holds for each $X_t^{\mathrm{ss}}$, $0 < t < 1$, and $X$ contains $\mathcal{L}\mathcal{K}^{-1}$-stable points.

"Upper semicontinuous" means: constant except at finitely many $t$, where $X_{t\pm\varepsilon}^{\mathrm{ss}} \subset X_t^{\mathrm{ss}}$. The $X_t^{\mathrm{ss}}$ will then satisfy a 'refinement lemma' akin to (1.2): each $X_t^{\mathrm{ss}} - X_{t\pm\varepsilon}^{\mathrm{ss}}$, if not empty, is a union of KN strata. (Indeed, if $X_{t\pm\varepsilon}^{\mathrm{ss}} \subseteq X_t^{\mathrm{ss}}$, then in fact $X_{t\pm\varepsilon}^{\mathrm{ss}}$ is the $(t\pm\varepsilon)$-semistable part of $X_t^{\mathrm{ss}}$, and application of Lemma (6.1) to $X_t^{\mathrm{ss}}$ reduces us to the projective case.) Parts (ii) and (iii) are automatic when $X$ is projective, but then, the conclusion is obvious, as noted earlier. We require, roughly, that no semistable sets "notice" the missing part of $X$. Note that semistability is delicate in the quasi-projective case ([MFK, Def. 1.7]).

*Proof.* Lemma (5.8) shows that the cohomology of $\mathcal{L}$ is the same for all $t < 1$. Because of (iv), $\mathcal{L}\mathcal{K}^{-1}$ is $G$-quasi-positive on $X_{1-\varepsilon}^{\mathrm{ss}}$, and the argument in (6.2) implies cohomology vanishing there. □

(6.6) *Example.* Let $X$ be the total space of $\mathcal{K}_\Sigma$, less the zero-section, on a smooth projective curve $\Sigma$ of genus 2. Since $X$ is quasi-affine, (6.5.i) is automatic. Further, (ii), (iii) and the first part of (iv) hold for $0 \le t < 1$, with $\mathcal{L} = \mathcal{K}_X$ and the natural $\mathbb{C}^\times$-action; but $H^1(\Sigma; \mathcal{K}) \neq 0$.

## 7. Differential forms and Hodge-to-de Rham spectral sequence

The two main results (7.1) and (7.3) were motivated mainly by the applications (9.9) and (9.10) to the moduli stack of principal bundles over a curve. Those will be further developed in a joint paper with S. Fishel and I. Grojnowski, leading to combinatorial applications.

*Statements.* Unless otherwise noted, $G$ is reductive and $X$ a completely KN-stratified $G$-manifold (§1). The first theorem involves the (strictly adapted) line bundle $\mathcal{L}$ and is the straightforward generalization of (2.11).

(7.1) THEOREM. *Restriction induces an isomorphism*

$$H^q\left(X_G; \Omega^p(\mathcal{L})\right) = H^q\left(X_G^{\mathrm{ss}}; \Omega^p(\mathcal{L})\right)$$

*of cohomologies over the smooth sites of the two stacks.*

(7.2) *Caution.* The "cohomology of $\Omega^p(\mathcal{L})$ over the smooth site of the stack" (cf. Remark 7.7) is the invariant hypercohomology of a complex $\mathcal{G}r^p_{\mathrm{Hodge}} \otimes \mathcal{L}$ over $X$; see the paragraph preceding (7.9). It is *not* the invariant part of $H^q(X; \Omega^p(\mathcal{L}))$, for which (7.1) can fail, as we saw for $\Omega^{\mathrm{top}}$. When $G$ acts freely on $X^{\mathrm{ss}}$, the second space in (7.1) equals $H^q(X^{\mathrm{ss}}/G; \Omega^p(\mathcal{L}))$, and, in this case, (7.1) generalizes (3.2.a). Already for finite stabilizers, we need to use the modules of *orbifold differentials* on $X^s/G$ [St, §1]. In general, the interpretation of $q_*^G \mathcal{G}r^p_{\mathrm{Hodge}}$ is more sophisticated than the definition: it should be the $p^{\mathrm{th}}$ Hodge-graded part of a *Hodge module* (in the sense of M. Saito), the total direct image of $\mathbb{C}$ along $q^G : X_G^{\mathrm{ss}} \to X^{\mathrm{ss}}/G$. Unfortunately, I know of no reference for equivariant mixed Hodge modules where this is made rigorous.

The more interesting case is $h = 0$, when the spaces in (7.1) constitute the $E_1$ term in the *Hodge-to-de Rham (Frölicher) spectral sequence for the stack $X_G$*. For classifying stacks of complex Lie groups, this was studied by Cathelineau [C], who viewed it as a holomorphic analogue of the Bott-Shulman spectral sequence. The following result (a) is quite clear when $X$ is proper, but less obvious in general.

(7.3) THEOREM. (a) *The Hodge-to-de Rham spectral sequence for $X_G$ collapses at $E_1$, which equals $\mathrm{Gr} H_G^*(X)$, for the algebro-geometric Hodge filtration. There is a Hodge decomposition*

$$H_G^n(X) \cong \bigoplus_{p+q=n} F^p_{\mathrm{Hodge}} H_G^n(X) \cap \overline{F^q_{\mathrm{Hodge}} H_G^n(X)}.$$

(b) *The same holds for the equivariant* HdR *sequences with supports on each* KN *stratum.*

(7.4) *Remark.* Kirwan [K1, Thm. 5.4] proved the *equivariant perfection* of the Morse function $\|\mu\|^2$, that is, the collapse at $E_1$ of the equivariant Cousin sequence of the stratification, with $\mathbb{C}$ coefficients. She further showed that the induced filtration on $H_G^*(X)$ was compatible with Hodge decomposition [K1, §15]:

$$\mathrm{Gr}^c H_G^{p,q}(X) = \bigoplus_{\mathrm{codim}\, S = c} H_L^{p-c, q-c}(Z^\circ).$$

We shall use these facts in proving (7.3) (see, however, Remark 7.22).

*The Cartan-Dolbeault complex.* Call $K$ a compact form of $G$, and let $\{A^\bullet(X), d\}$ be the $C^\infty$ de Rham complex of $X$. The *Cartan model* for $H_K^*(X)$ is the complex

$$A_K^\bullet(X) := \left[ A^\bullet(X) \otimes \mathrm{Sym}^\bullet \mathfrak{g}^t \right]^K$$

with differential

$$d_K = d + \sum_a \xi^a \cdot \iota(\xi_a) \ .$$

Here, $\xi_a$ is a basis of $\mathfrak{k}$, $\xi^a$ the dual basis of $\mathfrak{k}^t$, and $\iota$ denotes interior multiplication by the vector field generating the $\mathfrak{g}$-action. Elements of $\mathfrak{g}^t$ have degree 2. The complex structure refines the grading on $A^\bullet(X)$ to a bigrading, which we extend to $A_K^\bullet$ by declaring $\mathfrak{g}^t$ to have type $(1,1)$. The first (holomorphic) degree defines the descending *Hodge filtration* of $\{A_K^\bullet(X), d_K\}$:

(7.5) $$F_{\mathrm{Hodge}}^p A_K^\bullet(X) = \bigoplus_{r+s \geq p} \left[ A^{r,\bullet}(X) \otimes \mathrm{Sym}^s \mathfrak{g}^t \right]^K \ .$$

The associated graded complex $\bigoplus_p \mathrm{Gr}_{\mathrm{Hodge}}^p A_K^\bullet(X)$, defined by $r+s = p$ in (7.5), carries the *Cartan-Dolbeault differential* $\bar{\partial}_K := \bar{\partial} + \sum_a \xi^a \cdot \iota(\xi_a^{1,0})$. The spectral sequence resulting from (7.5) is the *Hodge-to-de Rham* (HdR) *spectral sequence for the stack $X_G$*.

For an invariant locally closed subvariety $S \subset X$, the HdR *spectral sequence for $X_G$ with supports in $S_G$* (converging to the equivariant cohomology of $X$ with supports in $S$) arises similarly, from the Hodge filtration on the total complex of the restriction arrow $A_K^\bullet(U) \to A_K^\bullet(U - S)$ when we use some invariant open set $U \subset X$ in which $S$ is closed.

(7.6) *Remark.* When $X$ is not proper, (7.5) computes the *analytic* HdR sequence, whereas it is the *algebraic* one that concerns us. One can show that the two agree, when $X$ is completely KN stratified,[7] but a simpler way around the difficulty is to choose a smooth $G$-compactification $\bar{X}$ of $X$, with normal-crossing complement, and replace $A^\bullet(X)$ by the Dolbeault resolution of the meromorphic extension to $\bar{X}$ of the holomorphic de Rham complex of $X$. (Taking the direct limit over all choices removes the dependence on $\bar{X}$.)

(7.7) *Remark.* The simplicial homotopy quotient variety $X_\bullet$, associated to the $G$-action by the bar construction, represents the stack $X_G$ in the simplicial theory of stacks. Its hypercohomology with constant coefficients is the $K$-equivariant cohomology of $X$. It follows from (7.3) that $\{A^\bullet(X_K), F^p, d_K\}$ is a

---

[7]One can show, in this case, that the algebraic sheaf cohomology is computed by the $K$-finite part of the Dolbeault complex.

*real, pure Hodge complex* [D, §6] representing $\mathbb{R}\Gamma(X_\bullet; \mathbb{C})$ in the derived category of (mixed, *a priori*) Hodge complexes. Note, when $X$ is not proper, that the filtration on $A^\bullet(X)$ by the holomorphic degree is *not* the Hodge filtration; so the result is not completely obvious.

(7.8) *Remark.* (7.3) implies the $\partial_K \bar{\partial}_K$-lemma for $A_K^\bullet$ [DGMS, 5.17]. The main theorem of that paper, to which we refer for definitions, implies that the differential graded algebra $\{A_K^\bullet(X), d_K\}$ is equivalent to its cohomology (with zero differential). Thus, subject to restrictions on $\pi_1$, the rational homotopy type of the stack $X_G$ (which is the homotopy quotient of $X$ by $G$) is determined by its cohomology $H_G^*(X)$, as are the rational homotopy types of inclusions $U_G \subset X_G$ of open unions of KN strata. In truth, the $\partial_K \bar{\partial}_K$-argument for open $U$ assumes that the analytical difficulty (7.6) has been properly addressed; however, the second proof of the main theorem in [DGMS], along with Morgan's follow-up work on open varieties, shows that purity of the Hodge structure and formality for complete $X$ suffice for the topological corollary.

*Algebraic description of $E_1$.* The bundle $\Omega_X^\bullet \otimes \mathrm{Sym}^\bullet \mathfrak{g}^t$ of Cartan-Kähler differentials, filtered by the Hodge degree, carries the equivariant de Rham operator $\partial_G = \partial + \sum_a \xi^a \cdot \iota(\xi_a)$. While $\partial_G^2 \neq 0$, the degree-zero part $\sum_a \xi^a \cdot \iota(\xi_a)$ is a differential, so that the hypercohomologies of the graded parts

$$\mathcal{G}r_{\mathrm{Hodge}}^p := \mathcal{G}r_{\mathrm{Hodge}}^p \left( \Omega_X^\bullet \otimes \mathrm{Sym}^\bullet \mathfrak{g}^t \right)$$

are defined. (Script $\mathcal{G}$ indicates sheaves, rather than vector spaces.)

(7.9) LEMMA. *The Cartan-Dobeault $E_1^{p,q}$ equals $\mathbb{H}_G^q \left( X; \mathcal{G}r_{\mathrm{Hodge}}^p \right)$, and $\partial_G$ gives rise to successive differentials which agree with those coming from* (7.5).

(7.10) *Remark.* One can show that $\mathbb{H}_G^q \left( X; \mathcal{G}r_{\mathrm{Hodge}}^p \right)$ agrees with the hypercohomology $\mathbb{H}^q(X_\bullet; \Omega^p)$ of $\Omega^p$ over the simplicial variety $X_\bullet$ of (7.7), and that the differentials arise from de Rham's operator.

*Proof.* $\mathbb{H}^*(X; \mathcal{G}r^p)$ is a complex-algebraic representation of $G$, so that

$$\mathbb{H}_G^* \left( X; \mathcal{G}r^p \right) = \mathbb{H}^* \left( X; \mathcal{G}r^p \right)^K;$$

the latter is computed by $(\mathrm{Gr}^p \Omega_K(X), \bar{\partial}_K)$, by exactness of $K$-invariants and by Dolbeault's theorem. $\square$

From this, one can see that the spectral sequence depends only on the quotient stack $X_G$, not on $X$ or $G$. In terms of equivariant cohomology, let $G' \supset G$ be reductive, and call $X'$ the induced space $G' \times^G X$. After factoring out the Koszul complex $\mathrm{Sym}^\bullet(\mathfrak{g}'/\mathfrak{g})^t \otimes \Lambda^\bullet(\mathfrak{g}'/\mathfrak{g})^t$, $\mathrm{Gr}\left(\Omega_{X'}^\bullet \otimes \mathrm{Sym}^\bullet \mathfrak{g}'^t\right)$ becomes quasi-isomorphic to the complex induced on $X'$ from $\mathcal{G}r\left(\Omega_X^\bullet \otimes \mathrm{Sym}^\bullet \mathfrak{g}^t\right)$.

(7.11) COROLLARY. *The HdR sequences for $X_G$ and $X'_{G'}$ are naturally isomorphic, from $E_1$ onwards.*

(7.12) *Remark.* Taking $G' = \mathrm{GL}_N$ yields an HdR spectral sequence for $X_G$ for any linear $G$; it can be shown to agree with "the" HdR sequence (7.10). The first part of (7.9) still holds; but when $G$ is *not* reductive, extra differentials arise from group cohomology (already when $X$ is a point; cf. [C]).

*Proof of* (7.3). Assume that $G$ is connected (the components only extract the invariant parts of the terms below). By (7.9), the $E_1$-term in the $G$-HdR sequence is

$$(7.13) \qquad E_1^{p,q} = \mathbb{H}^q\left(X; \mathcal{G}r_{\mathrm{Hodge}}^p\right)^G \Rightarrow H_G^{p+q}(X) .$$

For each $p$, filtering by the degree $r$ of the differential form gives yet another spectral sequence

$$(7.14) \qquad E_1^{r,s} = \left[H^s(X; \Omega^r) \otimes \left(\mathrm{Sym}^{p-r} \mathfrak{g}^t\right)\right]^G \Rightarrow \mathbb{H}^{r+s}\left(X; \mathcal{G}r_{\mathrm{Hodge}}^p\right) .$$

For proper $X$, collapse of the usual HdR implies that the $E_1$ term is $H^*(X) \otimes (\mathrm{Sym}^* \mathfrak{g}^t)^G$. This vector space already equals $H_G^*(X)$ (non-canonically), because the purity of the Hodge structure (or [K1, 5.8], in the projective case) forces the collapse of the Leray sequence for $X_G \to BG$ with complex coefficients. So there can be no further differentials in (7.13) or (7.14).

For any $X$, there is, for each Hodge degree $p$, a $G$-equivariant Cousin spectral sequence converging to the $E_1$ term in (7.13), involving the $G$-hypercohomologies with supports on the KN strata:

$$(7.15) \qquad E_1^{m,n} = \mathbb{H}_{S(m)}^{m+n}\left(X; \mathcal{G}r_{\mathrm{Hodge}}^p\right)^G \Rightarrow \mathbb{H}^{m+n}\left(X; \mathcal{G}r_{\mathrm{Hodge}}^p\right)^G ,$$

for a descending ordering $\{S(m)\}_{m \in \mathbb{N}}$ of the strata. Fix $m$ and let $q = m + n$; the left-hand side is the $E_q^{p,q}$ term of a $G$-HdR spectral sequence with supports on $S = S(c)$,

$$(7.16) \qquad E_1^{p,q} = \mathbb{H}_S^q\left(X; \mathcal{G}r_{\mathrm{Hodge}}^p\right)^G T \Rightarrow H_{S_G}^{p+q}(X_G; \mathbb{C})$$

abutting to equivariant cohomology, with supports on $S$, of the constant sheaf $\mathbb{C}$.

Collapse of all sequences (7.16) and the Hodge splitting in (7.3) will follow by induction on the rank of $G$. We shall see below how the inductive assumption forces the collapse of (7.16) for *unstable S*. Granting that, let us prove its collapse for $S = X^{\mathrm{ss}}$. Realize $X^{\mathrm{ss}}$ as the semistable stratum in some smooth projective $\bar{X}$, as in (6.1). With $\bar{X}$ replacing $X$, the $m \neq 0$ part of the Cousin sequence (7.15), and its abutment, are the $\mathrm{Gr}_{\mathrm{Hodge}}^p$ of their counterparts with $\mathbb{C}$ coefficients. But [K1, Thm. 5.4] asserts the collapse at $E_1$ of the latter sequence. As $\mathrm{Gr}_{\mathrm{Hodge}}^p$ is exact, the $m \neq 0$ part of (7.15) embeds into $E_\infty$. So there are no differentials, and the left edge $\mathbb{H}^n(X^{\mathrm{ss}}; \mathcal{G}r^p_{\mathrm{Hodge}})^G(X^{\mathrm{ss}}; \mathbb{C})$ agrees with $\mathrm{Gr}^p_{\mathrm{Hodge}} H^{n+p}_G(X^{\mathrm{ss}}; \mathbb{C})$. Hodge splitting for $X^{\mathrm{ss}}$ follows from the same for $\bar{X}$.

Collapse of (7.15) for any $X$ follows hence and from [K1, 5.4]: the differentials vanish in complex cohomology even prior to taking Hodge Gr's.

Let us explain the collapse of $G$-HdR for an unstable stratum $S$, mapping, as in Section 1, to $G \times^P Z^\circ$ under $\varphi$. As in Section 2, $\mathbb{H}^*_S\left(X; \mathcal{G}r^p_{\mathrm{Hodge}}\right)^G$ is the $G$-hypercohomology of a complex $\mathcal{R}_S \mathcal{G}r^p_{\mathrm{Hodge}}$ of sheaves of $S$. There is a natural inclusion

$$(7.17) \qquad \mathcal{G}r^p_{\mathrm{Hodge}}\left(\Omega^\bullet_S \otimes \mathrm{Sym}^\bullet \mathfrak{g}^t\right) \longrightarrow \mathcal{R}_S \mathcal{G}r^{p+c}_{\mathrm{Hodge}}\left(\Omega^\bullet_X \otimes \mathrm{Sym}^\bullet \mathfrak{g}^t\right)$$

($c = \mathrm{codim}(S)$), coming from the obvious maps (where $\mathcal{J}_S$ is the ideal sheaf of $S$),

$$(7.18) \ \Lambda^\bullet(T^t S) \to \Lambda^{\bullet+c}(T^t X) \otimes \Lambda^{\mathrm{top}}(T_S X), \quad \mathcal{O}/\mathcal{J}_S \otimes \Lambda^{\mathrm{top}}(T_S X) \to \mathcal{R}_S \mathcal{O}.$$

CLAIM A. *(7.17) induces a "Thom isomorphism" in $G$-hypercohomologies*,
$$(7.19)$$
$$\mathbb{H}^q\left(S; \mathcal{G}r^p_{\mathrm{Hodge}}\left(\Omega^\bullet_S \otimes \mathrm{Sym}^\bullet \mathfrak{g}^t\right)\right)^G \cong \mathbb{H}^{q+c}_S\left(X; \mathcal{G}r^{p+c}_{\mathrm{Hodge}}\left(\Omega^\bullet_X \mathrm{Sym}^\bullet \mathfrak{g}^t\right)\right)^G.$$

If so, it suffices to check the collapse of the $G$-HdR sequence for $S$ (rather than with supports in $S$). By (7.12), that is also the HdR sequence of the stack $Y^\circ_P$. We may assume the collapse at $E_1$ of HdR on $Z^\circ_L$ ($\mathbb{T}$ acts trivially on $Z^\circ$, so we are reduced to the lower-rank group $L/\mathbb{T}$); it thus suffices to prove:

CLAIM B. *The morphism of stacks $Y^\circ_P \to Z^\circ_L$ induced by $\varphi$ gives an isomorphism on $E_1$ terms.*

*Proof of the claims.* Freedom of the $\bar{\mathfrak{u}}$-action cancels a Koszul factor $\mathrm{Sym}^\bullet \bar{\mathfrak{u}}^t \otimes \Lambda^\bullet \bar{\mathfrak{u}}^t$ in $\varphi_* \mathcal{R}_S \mathcal{G}r^{p+c}_{\mathrm{Hodge}}$. Its "Gr" for the residue filtration is, then, the induced bundle, from $Z^\circ$ to $G \times^P Z^\circ$, of

$$(7.20) \quad \bigoplus_r \varphi_* \mathcal{G}r^{p+c-r}_{\mathrm{Hodge}}\left(\Omega^\bullet_Y \otimes \mathrm{Sym}^\bullet \mathfrak{p}^t\right) \otimes \mathrm{Sym}(T_S X) \otimes \Lambda^r(T_S X)^t \otimes \det(T_S X).$$

Combining the last two factors into $\Lambda^{c-r}(T_S X)$ shows that the $\mathbb{T}$-invariants are in $r = c$ and the constant line in $\mathrm{Sym}(T_S X)$. This picks the image of (7.17), proving (A). The surviving factor $\varphi_* \mathcal{G}r^p_{\mathrm{Hodge}}\left(\Omega^\bullet_Y \otimes \mathrm{Sym}^\bullet \mathfrak{p}^t\right)$ admits a further, $\mathbb{T}$-compatible degeneration to

$$(7.21) \quad \bigoplus_{r,s} \mathcal{G}r^{p+c-r-s}_{\mathrm{Hodge}}\left(\Omega^\bullet_Z \otimes \mathrm{Sym}^\bullet \mathfrak{l}^t\right) \otimes \mathrm{Sym}^r \mathfrak{u}^t \otimes \Lambda^s(T_Z Y)^t \otimes \mathrm{Sym}(T_Z Y)^t.$$

The $\mathbb{T}$-invariant part is the first factor, proving Claim B and completing the argument. □

*Proof of* (7.1). After twisting by the line bundle, the argument (7.16)–(7.19) implies the vanishing of all unstably supported cohomologies, because of the positive $\mathbb{T}$-weights on $\mathcal{L}$. □

(7.22) *Remark.* A more careful analysis shows that (7.3) consists of two portions:

- If $X$ is KN stratified, the Cousin sequence (7.15) collapses at $E_1$.

- If, additionally, (1.1.v) holds for a stratum $S$, the $G$-HdR spectral sequence with supports on $S$ collapses at $E_1$, and the filtration on $E_1 = E_\infty$ gives the Hodge splitting on cohomology.

Their proof is a refinement of Bott's argument [K1, 5.4] for the equivariant perfection of the stratification. However, being unaware of any application of the separate statements, I shall not prove them here.

## 8. Application to $G$-bundles over a curve

The final sections briefly explain how to recover a cohomology vanishing theorem of [T1], and prove vanishing of higher cohomologies of line bundles over moduli spaces of $G$-bundles over a curve. Also, Theorems 7.1 and 7.3 are restated for the stacks of bundles (Thms. 9.9 and 9.10), but serious applications will be discussed elsewhere. Everywhere, $G$ is reductive and connected.

The stack $\mathfrak{M}$ of algebraic $G$-bundles over a smooth projective curve $\Sigma$ of genus $g$ has several presentations. All global ones involve infinite-dimensional objects, which the equivariant cohomology of Appendix A cannot quite handle. Still, sheaf cohomology over $\mathfrak{M}$ (defined as cohomology over its smooth site) can be recovered as a limit over open substacks of finite type, which admit quasi-projective modulo reductive presentations. We can choose these substacks to be unions of finitely many *Shatz strata* (8.4).

One technical point deserves mention. The most explicit presentation, the Atiyah-Bott stack (8.3), is the *analytic stack underlying the stack of algebraic bundles.* To argue as in Section 2, we must use *algebraic* sheaf cohomology. That algebraic and analytic sheaf cohomologies agree is true but not obvious.[8] Avoiding that fact, the existence of algebraic presentations shows that algebraic sheaf cohomology is well-defined; while the analytic presentation (8.3) serves only to identify the subgroup $\mathbb{T}$ (of the gauge group) whose action on the strata and their normal bundles enforces the vanishing of local cohomologies.

---

[8]It follows from [T1], which reduced everything to Lie algebra cohomology. See Remark (7.6) for another argument.

(8.1) *Weil's adèle group presentation.* For semisimple $G$, $\mathfrak{M}$ is the double quotient stack $G(\Sigma)\backslash G_{\mathbb{A}}/G_{\mathbb{O}}$, where $\mathbb{A}$ is the adèle ring of $\Sigma$, $\mathbb{O}$ the subring of adèlic integers, and $G(\Sigma)$ the (algebraic ind-) group of $G$-valued rational maps. Alternatively, $\mathfrak{M}$ is the quotient under $G(\Sigma)$ of the ind-variety $G_{\mathbb{A}}/G_{\mathbb{O}}$, the *adèle flag variety.* Equality between the moduli and the double coset stack holds for any linear group at the level of $\mathbb{C}$-points; however, equality of the stacks requires Harder's result, and its refinement to families [DS], that algebraic bundles with *semisimple* structure group are trivial over an affine curve (this fails for $\mathrm{GL}_1$).[9]

(8.2) *The uniformization theorem* [LS]. One can economize in Weil's description by considering the affine curve $\Sigma^{\times}$ obtained from $\Sigma$ by removing a point $p$, and the groups $G[[z]]$ and $G((z))$ of formal holomorphic and of formal Laurent loops, using a formal coordinate centered at $p$. Then $\mathfrak{M}$ is equivalent to the quotient stack of the *local flag variety* $Q := G((z))/G[[z]]$ by the group $G[\Sigma^{\times}]$ of regular $G$-valued maps on $\Sigma^{\times}$. (Both are ind-objects.) Again, the key ingredient is [DS], so $G$ must be semisimple.

For another look at the same presentation, consider the quotient $X_{\Sigma^{\times}} := G((z))/G[\Sigma^{\times}]$. This is a scheme of infinite type, but locally it is fibered in (infinite-dimensional) affine spaces over a smooth finite-dimensional variety. (The fibers are normal subgroups of finite codimension in $G[[z]]$, and the base is a locally universal family of $G$-bundles over $\Sigma$, with suitable "level structure".) Then, $\mathfrak{M} \cong X_{\Sigma}/G[[z]]$. Similar statements hold for any number of punctures on $\Sigma$.

(8.3) *The Atiyah-Bott construction* [AB]. The identity component of $\mathfrak{M}$, in the *analytic* category, is the quotient stack of $\mathcal{A}$, the space of smooth, $\mathfrak{g}$-valued $(0,1)$-forms, under the gauge action of the group $\mathcal{G}(G)$ of smooth $G$-valued maps. (If $G$ is multiply connected, extra components of $\mathfrak{M}$ arise from connections and gauge transformations on various $C^{\infty}$ principal $G$-bundles over $\Sigma$.) Locally on $\mathcal{A}$, there are normal subgroups of finite codimension in $\mathcal{G}(G)$ acting freely, with smooth, finite-dimensional quotients; so this presentation is locally equivalent to a finite-dimensional one. It was suggested in [AB], and proved in [Da] and [R], that the Morse stratification of $\mathcal{A}$ by the Yang-Mills functional (defined after choice of a Kähler metric on $\Sigma$ and an inner product on $\mathfrak{g}$) agreed with the Shatz stratification (8.4).

In the spirit of Section 1, unstable strata are described as follows. A dominant 1-parameter subgroup $(1 - \mathrm{psg})$ $\gamma$ of $H \subset G$ determines $L$, $P$ and $\mathfrak{u}$ as in Section 1. Smooth $L$-bundles (equivalently, $P$-bundles) over $\Sigma$ are classified by $H^2(\Sigma; \pi_2(BL)) = H^2(\Sigma; \pi_1(L))$, so $\gamma^{-1}$ determines an equivalence class of smooth reductions to $P$ of our $G$-bundle. A Yang-Mills connection $A$ determines a $1-\mathrm{psg}$ $\mathbb{T}$ of $\mathcal{G}(G)$, whose infinitesimal generator is minus the curvature. Assume that $\mathbb{T}$ is locally conjugate to $\gamma$. (This happens if and only if $A$ determines an $L$-reduction of type $\gamma^{-1}$.) The $\mathbb{T}$-unstable stratum $S$ in $\mathcal{A}$ is the $\mathcal{G}(G)$-orbit of the set $Y^\circ$ of those connections which flow, as $t \to \infty$ in $\mathbb{T}$, into the set $Z^\circ$ of the *semistable* connections of topological type $\gamma^{-1} \in \pi_1(L)$. The connections in $Z^\circ$ are necessarily subordinate to the $L$-reduction of the $G$-bundle determined by $A$, those in $Y^\circ$ to the corresponding $P$-reduction. If we call $\mathcal{G}(P)$ and $\mathcal{G}(L)$ the obvious gauge groups, the quotient stack $\mathfrak{S} = S_{\mathcal{G}(G)}$ (a substack of $\mathfrak{M}$) is equivalent to $Y^\circ_{\mathcal{G}(P)}$, and classifies semistable $P$-bundles of type $\gamma^{-1}$; while $Z^\circ_{\mathcal{G}(L)}$ is the stack of semistable $L$-bundles of the same type.

---

[9]One could avoid [DS] by using nonstandard "configuration space" style ind-variety structures on $G(\Sigma)$ and $G_{\mathbb{A}}/G_{\mathbb{O}}$.

(8.4) *The Shatz stratification.* Shatz [Sh] defined functorial stratifications on algebraic families of vector bundles on $\Sigma$. Strata are indexed by concave polygonal functions with integral vertices, starting at $(0,0)$, the *Harder-Narasimhan* (HN) *polygons*. The vertex coordinates label the ranks and degrees of subbundles in a distinguished filtration with semistable subquotients. This filtration is a reduction of the structure group to a parabolic $P \subseteq \mathrm{GL}_n$; the HN polygon defines a homomorphism $\mathbb{T}: \mathbb{C}^\times \to \mathrm{GL}_n$, with generator corresponding to a dominant, regular weight $\beta$ of $P$, such that $\mathbb{T}^{-1}$ classifies the $P$-bundle topologically. The partial ordering on weights, defined by $\lambda > \mu$ if and only if $\lambda - \mu$ is in the span of negative roots, matches the partial ordering of the strata by the intersection relations between their closures [FM].

A similar description applies to any reductive $G$. Every unstable $G$-bundle has a unique parabolic reduction whose associated Levi bundle is semistable and for which the nilradical of the ad-bundle has positive HN slopes. The topological type of this reduction, described by a 1-parameter subgroup of $G$, labels the stratum to which the bundle belongs. The construction is functorial for group homomorphisms with finite kernels; for instance, the Shatz strata for $G$-bundles are pull-backs of those for $\mathrm{SL}(\mathfrak{g})$-bundles under the morphism $\mathfrak{M}_G \to \mathfrak{M}_{\mathrm{SL}(\mathfrak{g})}$ induced by the adjoint representation.

(8.5) *Line bundles over $\mathfrak{M}$.* For a torus $T$, $\mathfrak{M}_T \cong H^1(\Sigma; T) \times BT$, and Pic for each component is an abelian variety times $H^1(T; \mathbb{Z})$. For simple $G$, $\mathrm{Pic}(\mathfrak{M}) = H^2(\mathfrak{M}; \mathbb{Z})$, a finitely generated group of rank 1. In general, the exact sequence $[G, G] \to G \to G^{\mathrm{ab}} = T$ splits $\mathrm{Pic}(\mathfrak{M})$ into $\mathrm{Pic}(\mathfrak{M}_{[G,G]})$ and $\mathrm{Pic}(\mathfrak{M}_T)$, modulo finite groups (see (8.10) for more details). Call line bundles over $\mathfrak{M}_{[G,G]}$ *positive* if their lifts to the adèle flag variety are so. This is a convenient misnomer; there is no satisfactory notion of positivity for bundles over Artin stacks. Positivity over the $H^1(\Sigma; T)$-factor in $\mathfrak{M}_T$ has an obvious meaning. Positive line bundles over $\mathfrak{M}$ are those whose $\mathfrak{M}_{[G,G]}$ and $\mathfrak{M}_T$ components are positive. Similarly, one defines *semipositivity*.

Linearizing $\mathfrak{M}$ by any positive $\mathcal{L}$ leads to the same space

$$M = \operatorname{Proj} \bigoplus_{n \geq 0} \Gamma(\mathfrak{M}; \mathcal{L}^n) ,$$

the moduli space of semistable $G$-bundles. The morphism $q^{\operatorname{inv}} : \mathfrak{M}^{\operatorname{ss}} \to M$ from the open substack of semistable bundles to $M$ is known to be equivalent to a good quotient of a smooth quasi-projective variety by a reductive group (see [BR] and the references therein). $M$ is the quotient of $M_{[G,G]} \times M_{\operatorname{rad}(G)}$ by the finite group $H^1(\Sigma; \pi_1([G,G]))$, so, for cohomology vanishing questions, we only need to discuss simple groups.

In genus $g \geq 2$ (3, for $\mathfrak{sl}_2$) the complement of the stable points in $\mathfrak{M}$ has codimension 2 or more; (6.2) recovers the following result, due to Kumar and Narasimhan in the case when $\mathcal{L}$ descends to $M$.

(8.6) THEOREM (cf. [KN, Thm. B]). $H^*\left(M; q_*^{\operatorname{inv}}(\mathcal{K} \otimes \mathcal{L})\right) = 0$ *in positive degrees.*

(8.7) *Remarks.* (i) Kumar and Narasimhan deduce this directly from the theorem of Grauert and Riemenschneider; the key point is that $M$ is Gorenstein. (Cf. also Satz 2 in [Kn], following which $X//G$ is Gorenstein when (3.6) holds and $\mathcal{K}_G$ descends.)

(ii) (8.6) can be checked directly in genus 1, when the moduli space is a weighted projective space. The next section will show another way to remove the restriction on the genus.

The quantization theorem implies the following special case of [T1, Thm. 3]. The older statement about $H^0$ goes back to [BL] for $\operatorname{SL}_n$, and to [KNR] and [LS] in general.

(8.8) THEOREM. *For positive $\mathcal{L}$, $H^0(\mathfrak{M}; \mathcal{L}) = H^0(\mathfrak{M}^{\operatorname{ss}}; \mathcal{L})$, $H^{>0}(\mathfrak{M}; \mathcal{L})$ vanishes, as do all cohomologies with supports on unstable Shatz strata.*

*Proof.* We use the definitions and notation of (8.3); the objects essential to the argument (as opposed to those merely used to simplify the description) are algebraic. As we are only checking the signs of some $\mathbb{T}$-weights, an analytic description suffices.

At a $\mathbb{T}$-invariant connection $A \in \mathcal{A}$, defining a principal $L$-bundle $B_A$, the normal space $T_S \mathcal{A}$ (equal to $T_{\mathfrak{S}}\mathfrak{M}$, and replacing $T_S X$ in 2.4) is $H^1\left(\Sigma; B_A \times^L \bar{\mathfrak{u}}\right)$. Its $\alpha$-weight space has dimension $g-1+\langle\beta|\alpha\rangle$ and $\mathbb{T}$-weight $\langle\beta|\alpha\rangle$. The $\mathbb{T}$-weight on $\det(T_S \mathcal{A})$ is the sum of all $(g-1+\langle\beta|\alpha\rangle) \cdot \langle\beta|\alpha\rangle$. (It can happen, if $g=0$, that all dimensions $(g-1+\langle\beta|\alpha\rangle)$ vanish; $S$ is then open and there are no semistable bundles in the $A$-component of $\mathfrak{M}$. Otherwise, the weight on $\det(T_S \mathcal{A})$ is positive.)

Recalling that $Z^\circ_{\mathcal{G}(L)}$ is the $\mathbb{T}$-component of the stack of semistable $L$-bundles, let $E$ be the universal $L$-bundle over $\Sigma \times Z^\circ_{\mathcal{G}(L)}$ and $p$ the projection to $\Sigma$. As in Section 2, the cohomology with supports $H^*_{\mathfrak{S}}(\mathfrak{M}; \mathcal{V})$ of an algebraic vector bundle is computed from the cohomologies over $Z^\circ_{\mathcal{G}(L)}$ of $\mathcal{V} \otimes \det(T_S\mathcal{A}) \otimes \mathrm{Sym}(T_S\mathcal{A})$, tensored with the dual of the free commutative DGA on $\mathbb{R}p_*(E \times^L \mathfrak{u})[1]$; the first factor is the Gr of the residues, while the second is the fiberwise cohomology of the morphism $Y^\circ_{\mathcal{G}(P)} \to Z^\circ_{\mathcal{G}(L)}$ induced by the $\mathbb{T}$-flow. (The fibers are quotient stacks of $A^{0,1}(\Sigma; \mathfrak{u})$ by the $\bar{\partial}$-gauge action of $\mathcal{G}(U)$; this corresponds to the $U$-action on the fibers of $\varphi$ in Section 2.) All $\mathbb{T}$-weights on $\mathbb{R}p_*\left(U \times^L \mathfrak{u}\right)$ being negative, all cohomologies will vanish, if the $\mathcal{V}$-weights are nonnegative (or, even, not too negative). Since the $\mathcal{L}$-weight is positively proportional to $\|\beta\|^2$ (again, from the Atiyah-Bott description), we are done. $\square$

(8.9) *Remark.* Theorems 8.6 and 8.8 apply to all reductive groups. There seems to be no clean written account of the factorization formula (cf. 9.8.ii) for the dimension of $H^0$ for multiply connected groups, but an explicit "Verlinde formula" for all semisimple groups was recently derived by Alekseev, Meinrenken and Woodward [AMW].

(8.10) *Appendix: More on* $\mathrm{Pic}(\mathfrak{M})$ *for reductive* $G$. If $G$ is simply connected, $\mathfrak{M}$ is connected and $\mathrm{Pic}(\mathfrak{M})$ is free, of rank equal to the number of simple factors. For semisimple $G$, the components of $\mathfrak{M}$ are labeled by $\pi_1(G)$, and, for each component $\mathfrak{M}^{(p)}$, $p \in \pi_1(G)$, there appears in Pic the torsion subgroup of flat line bundles, dual to $\pi_1(\mathfrak{M}^{(p)}) = H^1(\Sigma; \pi_1(G))$, in addition to a free group of the expected rank. (Identification of the integer generators of the free group is a bit subtle; see [BLS], [T1].)

For any connected $G$, the group extension $G' = [G,G] \to G \to G^{\mathrm{ab}} = T$ leads, over a component $\mathfrak{M}^{(p)}_T$ of $\mathfrak{M}_T$ labelled by $p \in \pi_1(T)$, to a locally trivial fibration of stacks

$$\coprod_{p'} \mathfrak{M}^{(p')}_{[G,G]} \hookrightarrow \coprod_{p'} \mathfrak{M}^{(p')}_G \twoheadrightarrow \mathfrak{M}^{(p)}_T \ ,$$

with $p'$ ranging over the liftings of $p$ to $\pi_1(G)$. Here $\mathfrak{M}^{(p')}_G$ stands for the $p'$-component of $\mathfrak{M}_G$, but $\mathfrak{M}^{p'}_{[G,G]}$ is a component of a twisted version of the stack of $G'$-bundles, best described as the stack of "$\mathrm{Ad}(G')$-bundes with local $G'$-structures along $\Sigma$". The component in question is selected by the image of $p'$ in $\pi_1(\mathrm{Ad}(G'))$ (under the obvious map $G \to \mathrm{Ad}(G) = \mathrm{Ad}(G')$); and the obstruction to a global $G'$-structure is the further image in the center of $G'$. (The local $G'$-structures have the effect of cutting down the components of the gauge group, in the Atiyah-Bott presentation (8.3).) From the fibration of stacks, the Leray sequence with $\mathcal{O}^\times$ coefficients gives an exact sequence

$$1 \to \mathrm{Pic}\left(\mathfrak{M}_T^{(p)}\right) \to \mathrm{Pic}\left(\mathfrak{M}_G^{p'}\right) \to \mathrm{Pic}\left(\mathfrak{M}_{[G,G]}^{p'}\right) ,$$

which I claim is short and split.

We must only check this on the groups of components, because the (divisible) abelian variety in $\mathrm{Pic}(\mathfrak{M}_T)$ splits off. It suffices (by pure Hodge theory) to see that $H^2\left(\mathfrak{M}_G^{(p)}; \mathbb{Z}\right) \to H^2\left(\mathfrak{M}_{[G,G]}^{(p)}; \mathbb{Z}\right)$ is onto and that the torsion parts are isomorphic. Now, in the Leray sequence of the earlier fibration of stacks, with $\mathbb{Z}$ coefficients, the cokernel of this last map transgresses to the torsion-free group $H^3(\mathfrak{M}_T; \mathbb{Z})$, and rational splitting forces integral surjectivity. Equality of the torsion parts follows from $\pi_1(G') = \pi_1(G)^{\mathrm{tors}}$ and the earlier identification of $\pi_1(\mathfrak{M})$, which holds for any $G$.

## 9. Parabolic structures

There exist thorough treatments of moduli of vector bundles with parabolic structures ([MeS], [S2]), including the GIT variation problem ([BH], [Th]), but none, it seems, for parabolic $G$-bundles. This is partly because the definition of the moduli spaces adopted in the standard reference [BR] is "wrong" for certain parabolic weights (those on the far wall of the Weyl alcove, cf. (9.1); the true moduli space is then smaller than in [BR], and is always real-analytically isomorphic to a space of representations into the compact form of $G$). Most results extend from $\mathrm{GL}_n$ to other groups by embedding. For brevity assume that $G$ is simple.

A *quasi-parabolic datum* assigns to points $z_1, \ldots, z_m \in \Sigma$ parabolic[10] subgroups $\mathcal{P}_i$ of $G((z))$. The associated stack $\mathfrak{M}(\mathbf{z}, \mathbf{P})$ of $G$-bundles with quasi-parabolic structure is best described with respect to the adèlic presentation (8.2), as the quotient stack, under $G(\Sigma)$, of a modified adèlic flag variety, in which the local factor $G((z))/G[[z]]$ at $z_i$ is replaced by the *generalized flag variety* $G((z))/\mathcal{P}_i$. The *standard example* involves the subgroups $\mathcal{P}_i \subset G[[z]]$ of those formal holomorphic loops which take values, at $z = 0$, in specified parabolic subgroups $P_i \subset G$; $\mathfrak{M}(\mathbf{z}, \mathbf{P})$ is then the moduli stack of principal $G$-bundles over $\Sigma$, equipped with a reduction of the structure group to $P_i$ at $z_i$, and is an étale fiber bundle over $\mathfrak{M}$, with fiber $G/P_1 \times \cdots \times G/P_m$. In any case, $\mathfrak{M}(\mathbf{z}, \mathbf{P})$ is dominated by the stack $\mathfrak{M}(\mathbf{z}, \mathbf{B})$, defined by Borel data; more precisely, the latter is a fiber bundle over the former, with fiber a product of flag varieties of reductive subgroups of $G((z))$.

---

[10] *Parabolic subgroups* of $G((z))$ are those containing a conjugate of the *standard Borel subgroup*, the group of formal-holomorphic loops which become, when $z = 0$, a given Borel subgroup of $G$.

(9.1) *Remark.* When $G = \mathrm{SL}_n$, all generalized flag varieties are isomorphic to $G/P$-bundles over *twisted forms* of $G((z))/G[[z]]$, obtained by replacing $G((z))$ by multi-valued Laurent loops with period equal to a specified central element of $G$. Thus, inserting a twisted copy of $\mathrm{SL}_n((z))/\mathrm{SL}_n[[z]]$ leads to the stack of vector bundles with fixed, nontrivial determinant. For other groups and more general $\mathcal{P}_i$, the gauge condition at the $z_i$ which gives rise to $\mathfrak{M}(\mathbf{z}, \mathbf{P})$ is more awkward to spell out; and there is no morphism to the stack $\mathfrak{M}$ of $G$-bundles. These more general stacks must be included for uniform treatment (see Remark 9.3). This seems to have been missed in [BR].

Again, $\mathrm{Pic}(\mathfrak{M}(\mathbf{z}, \mathbf{P}))$ equals $H^2(\mathfrak{M}(\mathbf{z}, \mathbf{P})\mathbb{Z})$; we call line bundles over $\mathfrak{M}(\mathbf{z}, \mathbf{P})$ *(semi) positive* if they lift to (semi) positive ones over the modified adèlic flag variety. The class in $\mathbb{P}H^2(\mathfrak{M}(\mathbf{z}, \mathbf{P}); \mathbb{Q})$ of a positive line bundle is a *parabolic datum.*

(9.2) *Examples.* (i) In the standard example, $\mathrm{Pic}(\mathfrak{M}(\mathbf{z}, \mathbf{P}))$ splits into $\mathrm{Pic}(\mathfrak{M})$ and a full-rank subgroup of $\bigoplus_i \mathrm{Pic}(G/P_i)$, so that line bundles are determined by a *level* $k \in \mathbb{Z} \cong H^2(\mathfrak{M})$ and weights $\lambda_i$ of the parabolic subalgebras $\mathfrak{p}_i \subset \mathfrak{g}$. When $G$ is simply connected, all integral levels and weights are allowable.

(ii) Semi-positive bundles have $k \geq 0$ and dominant weights satisfying $\lambda_i \cdot \theta \leq k$ in the *basic* invariant inner product (in which the highest root $\theta$ has length $\sqrt{2}$). Positivity requires regular weights and strict inequalities. In general, positive line bundles require $\mathcal{P}_i$-dominant, regular affine weights of equal level.

(iii) The anticanonical bundle on $\mathfrak{M}(\mathbf{z}, \mathbf{P})$ is positive, in the sense just defined.

For positive $\mathcal{L}$, denote by $\mathfrak{M}(\mathbf{z}, \mathbf{P}, \mathcal{L})^{\mathrm{ss}} \subset \mathfrak{M}(\mathbf{z}, \mathbf{P})$ the complement of the base locus of large powers of $\mathcal{L}$. It depends on $\mathcal{L}$, but there are only finitely many of these semistable substacks ([BH] or [Th, §8] for $\mathrm{SL}_n$; argue by embedding in general). Now, $\mathfrak{M}(\mathbf{z}, \mathbf{P}, \mathcal{L})^{\mathrm{ss}}$ has a good quotient $q^{\mathrm{inv}} : \mathfrak{M}(\mathbf{z}, \mathbf{P}, \mathcal{L})^{\mathrm{ss}} \to M(\mathbf{z}, \mathbf{P}, \mathcal{L})$, the moduli space of bundles with parabolic structures;

$$M(\mathbf{z}, \mathbf{P}, \mathcal{L}) = \mathrm{Proj} \bigoplus_{n \geq 0} \Gamma\left(\mathfrak{M}(\mathbf{z}, \mathbf{P}); \mathcal{L}^n\right) .$$

(9.3) *Remark.* An inclusion of parabolic data $\mathbf{P} \subset \mathbf{P}'$ gives a map of flag varieties in the same direction, hence a morphism from $\mathfrak{M}(\mathbf{z}, \mathbf{P})$ to $\mathfrak{M}(\mathbf{z}, \mathbf{P}')$. This lifts positive line bundles to semipositive ones. Every semipositive line bundle over $\mathfrak{M}(\mathbf{z}, \mathbf{P})$ arises by lifting a positive one from a unique $\mathfrak{M}(\mathbf{z}, \mathbf{P}')$. For semipositive $\mathcal{L}$, we can define $M(\mathbf{z}, \mathbf{P}, \mathcal{L})$ as $M(\mathbf{z}, \mathbf{P}', \mathcal{L})$.

(9.4) *Caution.* If $g(\Sigma) > 0$, $M(\mathbf{z}, \mathbf{P}, \mathcal{L})$ is not empty precisely when $\mathcal{L}$ is semipositive (that is, semipositive = effective). In genus zero, the effective cone is more delicate to describe [TW]. We only need to know that $\mathcal{K}^{-1}$ is an *interior* point thereof, as soon as *three* Borel-marked points are present.

(9.5) *Stratifications of* $\mathfrak{M}(\mathbf{z}, \mathbf{P})$. (See [TW] for more details.) Call $\Sigma^\circ \subset \Sigma$ the affine complement of the marked points, $\hat{L}G$ the product of formal loop groups associated to them, $\hat{L}^+G \subset \hat{L}G$ the product of formal holomorphic loop groups, and $\mathbf{P}$ the product of the $\mathcal{P}_i$. If, as in (8.2), $X_{\Sigma^\circ} := \hat{L}G/G[\Sigma^\circ]$, then $\mathfrak{M}(\mathbf{z}, \mathbf{P}) \cong X_{\Sigma^\circ}/\mathbf{P}$.

Given a positive line bundle $\mathcal{L}$ of type $(k; \lambda_1, \ldots, \lambda_m)$, choose $N \in \mathbb{Z}^+$ so that each $N \cdot \lambda_i/k$ determines a homomorphism $\chi_i : \mathbb{C}^\times \to H$, using the basic inner product on $\mathfrak{g}$. Let $N\Sigma \to \Sigma$ be a cyclic covering of degree $N$, totally ramified at the $z_i$, call $\hat{L}_N G$ the $N\Sigma$-counterpart of $\hat{L}G$, and $\chi := (\chi_1, \ldots, \chi_m)$. Because $\mathbf{P} = \hat{L}G \cap \chi \cdot \hat{L}_N^+ G \cdot \chi^{-1}$, the embedding $E_\chi : X_{\Sigma^\circ} \subset X_{N\Sigma^\circ}$, defined by the inclusion $\hat{L}G \subset \hat{L}_N G$ followed by $\chi$-translation, gives a morphism from $\mathfrak{M}(\mathbf{z}; \mathbf{P})$ to the stack $\mathfrak{M}^N = X_{N\Sigma^\circ}/\hat{L}^+G$ of $G$-bundles over $N\Sigma$. The restriction there, under the Shatz strata of $\mathfrak{M}^N$, turns out to give a smooth KN stratification of $\mathfrak{M}(\mathbf{z}, \mathbf{P})$, to which $\mathcal{L}$ is strictly adapted.[11] Again, this is best seen in the Atiyah-Bott picture.

Call $\mathcal{A}$ the set of smooth $(0,1)$-connections on a $\Pi := \mathbb{Z}/N$-equivariant $G$-bundle over $N\Sigma$, which is trivial as a smooth bundle, but on which the $\Pi$-action translates the fiber over $z_i$ by $\exp(\lambda_i/k)$. If $\mathcal{G}$ is the smooth gauge group, $\mathfrak{M}(\mathbf{z}, \mathbf{P}) = \mathcal{A}^\Pi/\mathcal{G}^\Pi$. In a $\Pi$-invariant metric on $N\Sigma$, the entire Yang-Mills structure is $\Pi$-equivariant, and $\mathcal{A}^\Pi$ is smoothly stratified by the invariant parts of the Yang-Mills strata of $\mathcal{A}$. Each unstable $\mathcal{A}^\Pi$ stratum is the $\mathcal{G}^\Pi$-orbit of the set $S$ of $\Pi$-invariant connections, subordinate to a fixed $\Pi$-equivariant parabolic reduction, whose Levi parts define semistable, $\Pi$-equivariant $L$-bundles of fixed topological type. (Such a reduction is determined, after the choice of a $\Pi$-invariant Yang-Mills connection $A$ in $S$, by the curvature of $A$.) Note that $L$-bundles are topologically classified by $H^2_\Pi(\Sigma; C^\infty(L))$, so there is now a finite amount of information in addition to a 1-psg in $H$; so a stratum in $\mathcal{A}$ could break up into several strata in $\mathcal{A}^\Pi$. The normal space to a stratum within $\mathcal{A}^\Pi$, at some such $L$-bundle $B$, is the $\Pi$-invariant part of $H^1(\Sigma; B \times^L \bar{\mathfrak{u}})$, and the description of the local cohomologies proceeds as in the proof of (8.8).

*Remark.* A defter way to summarize the equivariant Atiyah-Bott construction, following [Bo] and [FS], involves the moduli stack of $G$-bundles over the orbifold $N\Sigma/\Pi$. Its components are labeled by assigning, to each marked point $z_i$, a conjugacy class of $N^{\text{th}}$ roots of unity in $G$. Corresponding to $\{\exp(\lambda_i/k)\}$ is our stack $\mathfrak{M}(\mathbf{z}, \mathbf{P})$. The basic statement (proved by this reduction to the equivariant picture) is that "Atiyah-Bott applies to compact Riemann orbi-surfaces."

---

[11] Note that the fundamental line bundle pulls back to $\mathcal{L}^{N/k}$. This functoriality is a general feature of Kirwan's stratification.

(9.6) THEOREM. *For semipositive $\mathcal{L}$,*

$$H^q(\mathfrak{M}(\mathbf{z}, \mathbf{P}); \mathcal{L}) = H^q\left(M(\mathbf{z}, \mathbf{P}, \mathcal{L}); q_*^{\text{inv}} \mathcal{L}\right),$$

*and vanishes if $q > 0$.*

(9.7) *Remark.* When $G = \text{SL}_2$, Mehta and Ramadas [MR] describe the cohomology over the moduli spaces in characteristics $p \neq 2, 3$. (If $p \neq 0$, the curve $\Sigma$ must be suitably generic.)

*Proof.* $\mathcal{L}$ is positive on some $\mathfrak{M}(\mathbf{z}, \mathbf{P}')$, and equality of the cohomologies follows, as in (8.8), from the properties of the stratification in (9.5). For the vanishing, we argue as in (6.5), after adding enough Borel data. The general quotients will then contain stable points, and (3.6) is automatic, because all semistable isotropies are abelian. Lemma 5.8 implies that

$$H^q\left(M(\mathbf{z}, \mathbf{B}, \mathcal{L}); \mathcal{L}\right) = H^q\left(M(\mathbf{z}, \mathbf{B}, \mathcal{L}\mathcal{K}^{-1}); \mathcal{L}\right).$$

(Note, in genus 0, that if $\mathcal{L}$ is effective and we add three or more Borel markings, then $\mathcal{L}\mathcal{K}^{-1}$ lies inside the effective cone; while, if $\mathcal{L}$ is not effective, $M \neq \emptyset$ and there is nothing to prove.) Finally, vanishing of higher cohomology over $M(\mathbf{z}, \mathbf{P}, \mathcal{L}\mathcal{K}^{-1})$ follows from (6.4). $\square$

(9.8) *Remarks.* (i) One recovers the dominant case of the cohomology vanishing in [T1, Thm. 3] by pushing semipositive line bundles over $\mathfrak{M}(\mathbf{z}, \mathbf{B})$ down to $\mathfrak{M}$; the result is a twist of a positive line bundle by a product of *evaluation vector bundles.*

(ii) Unlike the proof in [T1], this argument does not determine the dimension of the space of sections. However, (9.6) in genus zero implies Verlinde's factorization formula and determines the fusion rules; for it recovers the Lie algebra cohomology vanishing result on which [T1] is based (see [T2] for a short survey). One can of course also invoke the original factorization proof for $H^0(\mathfrak{M}(\mathbf{z}, \mathbf{P}))$ [TUY]. Note, finally, that [MW] recovers, by symplectic methods, the factorization formula for the *index* of $\mathcal{L}$ over (shift desingularizations of) $M(\mathbf{z}, \mathbf{P}; \mathcal{L})$. Direct computations of this index have also been given.

(iii) The proofs in [T1] relied on the "Borel-Weil-Bott" theorem of Kumar and Mathieu. The vanishing theorem for the genus zero stack, with parabolic structures at *two* marked points, *recovers* that theorem as a consequence of the BGG resolution. Of course, in finite dimensions, it is well-known that BGG implies the geometric Borel-Weil-Bott theorem; but the argument uses the Peter-Weyl theorem for $G$, which fails for $G((z))$. Use of the stack circumvents that difficulty.

(9.9) THEOREM. *The Hodge-to-de Rham spectral sequences for $\mathfrak{M}$, for $\mathfrak{M}(\mathbf{z}, \mathbf{P})$, or for any open substacks which are unions of Shatz strata, collapse at $E_1$ and induce the Hodge filtration on cohomology.*

The dimensions of $H^q(\mathfrak{M}; \Omega^p)$ thus equal the Hodge numbers of $\mathfrak{M}$; for those, see [T1, §5]. The Hodge structure of $\mathfrak{M}(\mathbf{z}, \mathbf{P})$ can be found by the same methods.

(9.10) THEOREM. *For positive $\mathcal{L}$,*

$$H^*\left(\mathfrak{M}; \Omega^p(\mathcal{L})\right) = H^*\left(\mathfrak{M}^{\mathrm{ss}}; \Omega^p(\mathcal{L})\right) ,$$

*and similarly for $\mathfrak{M}(\mathbf{z}, \mathbf{P})$.*

The spaces $H^q\left(\mathfrak{M}; \Omega^p(\mathcal{L})\right)$ will be discussed elsewhere. For $p = 0$, they are the famous *conformal blocks*, with dimensions given by Verlinde's factorization rule. At another extreme, in genus zero, $\mathfrak{M}^{\mathrm{ss}} = BG$, to which $\mathcal{L}$ restricts trivially; so $H^q\left(\mathfrak{M}; \Omega^p(\mathcal{L})\right) = H^{p,q}(BG)$. Note, in comparing with $\mathcal{L} = \mathcal{O}$, that the factor $H^*(\Omega G)$ has been lost from

$$H^*(\mathfrak{M}; \Omega^*) = H^*(\mathfrak{M}; \mathbb{C}) = H^*(BG) \otimes H^*(\Omega G) .$$

## 10. Complements in positive characteristic

Appropriately stated, the main theorems and some of their consequences in Sections 3 and 5 carry over to positive characteristic. The key point is that the description in Section 1 of the Kirwan-Ness stratification, and its properties (i)–(v), apply to smooth projective varieties over a perfect ground field $k$ [K1, §12].

One subtlety, when $\operatorname{char}(k) = p > 0$, concerns the distinction between *linearly reductive* and *geometrically reductive* groups. For the former, the category of representations is semisimple; the latter verify the weaker condition that any surjection $V \to k$ of $G$-modules splits after going over to a symmetric power. (These, plus affinity of $G$, may be taken as definitions.) It is known that $G$ is geometrically reductive if and only if its identity component $G_0$ has trivial unipotent radical, linearly reductive if $G_0$ is a torus and the order of the group of components is not a multiple of $p$ (see [MFK] for references).

Another complication is the existence of nonreduced group schemes; they appear naturally as stabilizers (for instance, the group scheme of $p^{\mathrm{th}}$ roots of unity is the kernel of the $p^{\mathrm{th}}$ power automorphism of $\mathbb{G}_m$). Luckily, they can be treated on a par with the groups, with geometric and linear reductivity defined as above. An important fact is that closed subgroup schemes of linearly reductive groups are also linearly reductive, while closed subgroup schemes of geometrically reductive ones are geometrically reductive if and only if the quotient space is affine. (The argument in [H] works for group schemes.) Note, though, that finite quotient stacks are not of Deligne-Mumford type, if the isotropies are not reduced.

GIT produces good quotients (3.1) $q : X^{ss} \to X^{ss}/G$ for geometrically reductive groups. Moreover, $q$ separates invariant closed sets, and each fiber contains a unique closed orbit, whose isotropy is a geometrically reductive group scheme. However, $q_*^G$ is only left exact; it is exact, if the following condition holds:

(L) Closed orbits in $X^{ss}$ have *linearly* reductive isotropy group schemes.

Luna's slice theorem then applies. The quantization theorem runs as follows.

(10.1) THEOREM. *Assumptions are as in* 3.2.a.
(a) $H_G^*(X; \mathcal{L}) = H_G^*(X^{ss}; \mathcal{L})$.
(b) *If* (L) *holds,* $H_G^*(X; \mathcal{L}) = H^*(X//G; \mathcal{L})$.
(c) *If $G$ is linearly reductive,* $H^*(X; \mathcal{L})^G = H^*(X//G; \mathcal{L})$.

The argument in Section 2 carries over, with the following changes:

- For the Gr of the sheaf of residues, $\mathrm{Sym}(T_S X)$ is replaced by the *divided symmetric algebra*;

- Lie algebra cohomology of $\mathfrak{u}$ is replaced by group cohomology of $U$ (resolved by the bar complex);

- $L$-invariants are replaced by $L$-cohomologies; their vanishing follows from the identity $H_L^*(V) = H_{L/\mathbb{T}}^*(V^{\mathbb{T}})$ (a consequence of the Hochschild-Serre spectral sequence) and from the vanishing of $\mathbb{T}$-invariants.

Most useful consequences require condition (L). Without it, the relative quantization theorem (5.2) takes the less attractive form $\mathbb{R}p_* \mathbb{R}q_*^{\prime G} \mathcal{M}^h = \mathbb{R}q_*^G \mathbb{R}\pi_* \mathcal{M}^h$. Subject to (L), it holds as originally stated, and Theorems 5.4–5.6 apply. Boutot's theorem (5.7) takes the following form:

(10.2) PROPOSITION. *If $X$ is smooth and* (L) *holds, there exists a smooth $\tilde{X}$, projective over $X$, such that $\tilde{X}//G$ has finite (linearly reductive) quotient singularities, $p : \tilde{X}//G \to X//G$ is birational and $\mathbb{R}p_* \mathcal{O} = \mathcal{O}$.*

By Grothendieck duality, $\mathbb{R}p_* \mathbb{D}_{\tilde{X}//G} \mathcal{O} = \mathbb{D}_{X//G} \mathcal{O}$ on dualizing complexes; so $p_* \tilde{\omega} = \omega$ on dimensional grounds, in particular, $X^{ss}/G$ is Cohen-Macaulay. This is the *Hochster-Roberts theorem* (see e.g. [Ke2]).

## Appendix. Quotient stacks

The theory of coherent sheaves over 1-stacks of finite type is described in [LM]. However, the case of quotients of varieties by group actions reduces to the equivariant cohomology of coherent sheaves, so similar to its topological counterpart that an economical account is possible (at least over $\mathbb{C}$).

(A.1) *Equivariant sheaves.* Given a complex scheme $X$ of finite type, acted upon by a linear algebraic group $G$, a (*quasi*) *coherent sheaf* over the quotient stack $X_G$ is a $G$-equivariant (quasi) coherent sheaf over $X$. Algebraicity of the group action can be imposed via a *cocycle condition* [MFK, Def. 1.6], or geometrically, by requiring that the induced action on the associated linear space Spec Sym of the sheaf be algebraic.

(A.2) *Categorical properties.* Equivariant quasi-coherent sheaves over $X$ form an abelian category $\mathfrak{Qcoh}^G(X)$ (or $\mathfrak{Qcoh}(X_G)$) under $G$-morphisms.[12] When $X$ is a principal $G$-bundle, $\mathfrak{Qcoh}^G(X)$ is equivalent to the category of sheaves over its base. When $X$ is a point $*$, quasi-coherent sheaves over $*_G$ (also denoted $BG$) are locally-finite algebraic representations of $G$. The following are easily checked. Call $p$ the projection $G \times X \to X$. For any $\mathcal{F} \in \mathfrak{Qcoh}(X)$, $p_* p^* F$ has a natural $G$-structure. The adjunction morphism $\mathcal{F} \to p_* p^* \mathcal{F}$ is equivariant, if $\mathcal{F}$ carries a $G$-action to begin with. More precisely, $p_* p^*$ is right-adjoint to the forgetful functor from $\mathfrak{Qcoh}(X_G)$ to $\mathfrak{Qcoh}(X)$. It therefore takes injectives to injectives, and it follows that $\mathfrak{Qcoh}^G(X)$ has enough injectives.

(A.3) *Cohomology.* Global sections of $\mathcal{F}$ over $X$ form an algebraic, locally-finite representation of $G$. The functor $\mathcal{F} \mapsto \Gamma(X_G; \mathcal{F}) := \Gamma(X; \mathcal{F})^G$ of invariant global sections is left exact, and its right derived functors are the *cohomology groups* of $\mathcal{F}$ over $X_G$, or the *equivariant sheaf cohomologies over $X$*. Over $BG$, $\Gamma^G$ selects the invariant vectors in a representation. Its derived functors are the *group cohomologies*, denoted $H^*(BG; \_)$ or $H^*_G(\_)$. Grothendieck's spectral sequence for the composition of functors

$$(A.4) \qquad E_2^{p,q} = H_G^p(H^q(X; \mathcal{F})) \Rightarrow H^{p+q}(X/G; \mathcal{F})$$

can be viewed as Leray's spectral sequence for the morphism $X_G \twoheadrightarrow BG$.

(A.5) *Shapiro's lemma.* Given an embedding $G \subset G'$, form the *induced space* $X' := G' \times^G X$. Restriction to $X$ gives an equivalence between the categories of (quasi) coherent equivariant sheaves over $X'_{G'}$ and $X_G$, commuting with the functor of invariant global sections; so sheaf cohomologies over the two stacks agree. (More generally, for a $G$-equivariant principal $G'$-bundle $P$ over $X$, the stacks $P_{G \times G'}$ and $X_G$ are equivalent, in the same sense.) When $X$ is a point, this is a classical result of Shapiro's: the $G'$-cohomology of a representation induced from $G$ agrees with the $G$-cohomology of the original. A reductive group $G'$ has no higher group cohomology, so $H^*_G(X; \mathcal{F}) = H^*(X'; \mathcal{F}')^{G'}$, in obvious notation. Reduction to invariants is a shortcut for extending standard theorems to the equivariant setting. Most important for us is the equivariant Cousin-Grothendieck spectral sequence,

$$(A.6) \qquad E_1^{c,d} = H^{c+d}_{S(c)_G}(X_G; \mathcal{F}) \Rightarrow H^{c+d}(X_G; \mathcal{F}) \,,$$

where $S(c)$ is a decreasing $G$-stratification of $X$ by locally closed subsets; the formula reads "cohomology of $\mathcal{F}$ over the stack $X_G$ with supports on the locally closed substack $S(c)_G$", or better, "equivariant cohomology of $\mathcal{F}$ over $X$ with supports in $S(c)$". It can be defined, after a reductive embedding of $G$, as invariant cohomology with supports, on the induced stratified space; alternatively, when $S(c)$ is closed in some open $U \subset X$, with open complement $j : V \hookrightarrow U$, it is the $G$-hypercohomology of the complex $\mathcal{F} \to \mathbb{R}j_* j^* \mathcal{F}$ over $U$ ("$G$-cohomology over $U$ relative to $V$").

---

[12] "Smooth descent of quasi-coherent sheaves" asserts that this category is equivalent to that of quasi-coherent sheaves over the *smooth site of the quotient stack $X/G$*. See [LM] for more details along these lines.

(A.7) *Étale slices.* Let $G$ be reductive, $X$ affine and normal, with quotient $q : X \to X/G$. Each fiber of $q$ contains a single closed orbit, in which the stabilizer $R$ of any point $x$ is reductive. Luna's theorem [MFK, App. A] asserts the existence of an *étale slice* $W$ through $x$, an affine $R$-subvariety of $X$ for which the morphism $G \times_R W \to K$ is *strongly étale*. When $X$ is smooth, $W$ is $R$-isomorphic to an étale neighborhood of $0$ in the normal bundle to $Gx$. "Strongly étale" means that the morphisms $q^G : X_G \to X/G$ and $q^R : W_R \to W/R$ on quotient stacks fit into an obvious Cartesian square, in which (étale-locally near $x$) the quotient spaces are isomorphic and the quotient stacks are equivalent. Note that $\mathcal{K}_G$ on $X$ corresponds to $\mathcal{K}_R$ on $W$.

(A.8) *Deligne-Mumford stacks.* Morally, DM stacks arise by gluing together finite quotient stacks via Shapiro equivalences.[13] Sheaves are obtained by compatible gluings of equivariant sheaves. Properties (such as coherence or smoothness) which can be tested locally have obvious meanings; so does the Kähler condition on a (globally defined) 2-form. Thus, when $X$ is a compact Kähler manifold with Hamiltonian $G$-action, a holomorphic slice argument [Sj1] shows that $X_G^s$ is a smooth DM stack, and carries a Kähler form obtained by symplectic reduction.

The (awkward) compatibility condition on Shapiro gluings (and the pathologies of nonseparated stacks) are best circumvented in Artin's presentation of DM stacks by *étale groupoids* $X_1 \rightrightarrows X_0$, groupoids in the category of algebraic spaces, where the two structural projections are étale morphisms. (One may allow smooth structure maps, provided the *diagonal morphism* $X_1 \to X_0 \times X_0$ is unramified: a slice argument produces an étale presentation.) Now, $X_G$ corresponds to the groupoid $G \times X \rightrightarrows X$, and is a separated DM stack if and only if the action is proper. A notion of equivalence of groupoids, extending (A.5), generates an equivalence for stacks; this agrees with the equivalence defined category-theoretically in [LM]. A stack is *separated* if the diagonal morphism is finite; it then has a separated algebraic space quotient [KM]. It is *proper* if it is separated and its quotient space is proper. Sheaves on the stack can be pushed down to the quotient space; the higher direct images vanish when the orders of the isotropies are invertible in the stalks.

---

[13]This works for separated stacks; the correct general definition starts from *groupoids*, as in the next paragraph.

For smooth stacks, we can define local and global $C^\infty$ functions, de Rham and Dolbeault complexes. The cohomology of (holomorphic bundles over) a smooth, separated DM stack can be computed from the global sections of de Rham's (Dolbeault's) complex, whose push-down to the quotient space is a complex of soft resolution sheaves. A Kähler structure determines Hilbert spaces of $L^2$ forms, and the differential-geometric Kähler package applies. In the smooth, compact, Kähler case, familiar theorems on manifolds (representation of cohomology by harmonic forms, the classical vanishing theorems of Kodaira-Nakano and Grauert-Riemenschneider and the collapse at $E_1$ of the Hodge spectral sequence) follow in the usual manner. This has been known for some time, albeit in the language of "$V$-manifolds" (see for instance [St, §1]); nowadays it is often attributed to "physicists".

(A.9) *Cohomology with (formally) proper supports* (see Deligne's appendix to [Ha]). Completing $X$ to a proper variety $X'$, any coherent sheaf $\mathcal{F}$ on $X$ admits a coherent extension $\mathcal{F}'$ over $X'$. Call $\hat{\mathcal{F}}'_Z$ the completion of $\mathcal{F}'$ along $Z := X' - X$. The cohomology $H^*_c(X; \mathcal{F})$ of $\mathcal{F}$ over $X$ with formally proper supports is the hypercohomology of the complex $\mathcal{F} \to \hat{\mathcal{F}}'_Z$ over $X'$. It depends only on $X$ and and $\mathcal{F}$; indeed, Serre duality asserts that $H^*_c(X; \mathcal{F})$ is the algebraic dual of $\mathrm{Ext}_X^{\dim X - *}(\mathcal{F}; \mathbb{D}\mathcal{O})$, where $\mathbb{D}\mathcal{O}$ is the dualizing complex of $X$. If $X$ and $\mathcal{F}$ carry a $G$-action, we can choose $X'$ and $\mathcal{F}'$ to do likewise; $H^*_c(X; \mathcal{F})$ is then the full dual of a locally finite $G$-representation, and Serre duality is equivariant.

More generally, one defines, for a compactifiable $G$-morphism $f: X \to Y$, an equivariant direct image with proper supports $\mathbb{R}f_!$, whose value is a pro-object in $D^b\mathfrak{Coh}^G(Y)$, the bounded derived category. For a good quotient $q: X^\circ \to X^\circ/G$, let $\mathbb{R}q_!^G$ be the downward dim $G$-shifted invariant part of $\mathbb{R}q_!$. Relative duality $\mathcal{E}\mathrm{xt}_{X^\circ/G}(\mathbb{R}q_!^G \mathcal{F}; \omega) = q_*^G \mathcal{E}\mathrm{xt}_{X^\circ}(\mathcal{F}; \mathcal{K}_G)$ shows that $\mathbb{R}q_!^G$ lands, in fact, in $D^b\mathfrak{Coh}(X^\circ/G)$.

STANFORD UNIVERSITY, STANFORD, CA
*E-mail address*: teleman@math.stanford.edu